\documentclass[12pt]{article}
\usepackage{graphicx}
\usepackage{amsmath,amsthm,amssymb,enumerate}%, esint}
\usepackage{euscript,mathrsfs}
\usepackage{color}
\usepackage{dsfont}
\usepackage{url}
\usepackage[left=2cm,right=2cm,top=3.5cm,bottom=3.5cm]{geometry}
\usepackage{color}
\usepackage[framemethod=tikz]{mdframed}
\allowdisplaybreaks
\usepackage{esint}

\newcommand{\avintO}[1]{\fint_{\Omega} #1 \dx}

\usepackage{soul}

\catcode`\@=11 \@addtoreset{equation}{section}

\catcode`\@=12

\newtheorem{Theorem}{Theorem}[section]
\newtheorem{Proposition}[Theorem]{Proposition}
\newtheorem{Lemma}[Theorem]{Lemma}
\newtheorem{Corollary}[Theorem]{Corollary}

\theoremstyle{definition}
\newtheorem{Definition}[Theorem]{Definition}

\newtheorem{Remark}[Theorem]{Remark}

\newcommand{\bTheorem}[1]{
	\begin{Theorem} \label{T#1} }
	\newcommand{\eT}{\end{Theorem}}

\newcommand{\bProposition}[1]{
	\begin{Proposition} \label{P#1}}
	\newcommand{\eP}{\end{Proposition}}

\newcommand{\bLemma}[1]{
	\begin{Lemma} \label{L#1} }
	\newcommand{\eL}{\end{Lemma}}

\newcommand{\bCorollary}[1]{
	\begin{Corollary} \label{C#1} }
	\newcommand{\eC}{\end{Corollary}}

\newcommand{\bRemark}[1]{
	\begin{Remark} \label{R#1} }
	\newcommand{\eR}{\end{Remark}}

\newcommand{\bDefinition}[1]{
	\begin{Definition} \label{D#1} }
	\newcommand{\eD}{\end{Definition}}

\newcommand{\Del}{\Delta_x}

\newcommand{\Ds}{\mathbb{D}_x}

\newcommand{\MTC}{\mathcal{T}}

\newcommand{\bfphi}{\boldsymbol{\varphi}}

\newcommand{\bFormula}[1]{
	\begin{equation} \label{#1}}
	\newcommand{\eF}{\end{equation}}

\newcommand{\vuh}{\vu_h}

\newcommand{\Divh}{{\rm div}_h}
\newcommand{\Gradh}{\nabla_h}

\newcommand{\Ov}[1]{\overline{#1}}

\newcommand{\Curl}{{\bf curl}_x}

\newcommand{\vr}{\varrho}
\newcommand{\vre}{\vr_\ep}

\newcommand{\vte}{\vt_\ep}
\newcommand{\vue}{\vu_\ep}
\newcommand{\tvr}{\wtilde \vr}
\newcommand{\tvu}{{\wtilde \vu}}
\newcommand{\tvt}{\wtilde \vt}

\newcommand{\vt}{\vartheta}
\newcommand{\vu}{\vc{u}}

\newcommand{\vc}[1]{{\bf #1}}

\newcommand{\Div}{{\rm div}_x}
\newcommand{\Grad}{\nabla_x}

\newcommand{\dx}{\,{\rm d} {x}}

\newcommand{\dt}{\,{\rm d} t }

\newcommand{\vU}{\vc{U}}

\newcommand{\intO}[1]{\int_{\Omega} #1 \ \dx}

\newcommand{\D}{{\rm d}}

\newcommand{\ep}{\varepsilon}

\newcommand{\R}{\mathbb{R}}
\newcommand{\vtB}{\vt_B}

\newcommand{\br}{ \nonumber \\ }

\def\softd{{\leavevmode\setbox1=\hbox{d}%
		\hbox to 1.05\wd1{d\kern-0.4ex{\char039}\hss}}}
\definecolor{Cgrey}{rgb}{0.85,0.85,0.85}
\definecolor{Cblue}{rgb}{0.50,0.85,0.85}
\definecolor{Cred}{rgb}{1,0,0}
\definecolor{fancy}{rgb}{0.10,0.85,0.10}
\definecolor{amaranth}{rgb}{0.9, 0.17, 0.31}

\newcommand\Cbox[2]{%
	\newbox\contentbox%
	\newbox\bkgdbox%
	\setbox\contentbox\hbox to \hsize{%
		\vtop{
			\kern\columnsep
			\hbox to \hsize{%
				\kern\columnsep%
				\advance\hsize by -2\columnsep%
				\setlength{\textwidth}{\hsize}%
				\vbox{
					\parskip=\baselineskip
					\parindent=0bp
					#2
				}%
				\kern\columnsep%
			}%
			\kern\columnsep%
		}%
	}%
	\setbox\bkgdbox\vbox{
		\color{#1}
		\hrule width  \wd\contentbox %
		height \ht\contentbox %
		depth  \dp\contentbox
		\color{black}
	}%
	\wd\bkgdbox=0bp%
	\vbox{\hbox to \hsize{\box\bkgdbox\box\contentbox}}%
	\vskip\baselineskip%
}

\mdfdefinestyle{MyFrame}{%
	linecolor=black,
	outerlinewidth=1pt,
	roundcorner=5pt,
	innertopmargin=\baselineskip,
	innerbottommargin=\baselineskip,
	innerrightmargin=10pt,
	innerleftmargin=10pt,
	backgroundcolor=white!20!white}

\newcommand{\tbf}{\textbf}

\newcommand{\mbf}{\mathbf}

\newcommand{\mfr}{\mathfrak}

\newcommand{\veps}{\varepsilon}

\newcommand{\wtilde}{\widetilde}

\renewcommand{\o}{\omega}

\newcommand{\lan}{\langle}
\newcommand{\ran}{\rangle}

\newcommand{\aleq}{\lesssim}

\renewcommand{\bot}{{\rm bot}}
\renewcommand{\top}{{\rm top}}

\allowdisplaybreaks

\begin{document}

%%%%%%%%%%%%%%%%%%%%%%%%%%%%%%%%

\title{\bf On the Rayleigh--B\' enard convection problem for rotating fluids}

\author{Francesco Fanelli 
\thanks{The work of F.F. has been partially supported by the project CRISIS (ANR-20-CE40-0020-01), operated by the French National Research Agency (ANR),
by the Basque Government through the BERC 2022-2025 program and by the Spanish State Research Agency through the BCAM Severo Ochoa excellence accreditation
CEX2021-001142. The author also aknowledges the support of the European Union through the COFUND program [HORIZON-MSCA-2022-COFUND-101126600-SmartBRAIN3].}
\and
Eduard Feireisl \thanks{The work of E.F. was partially supported by the
		Czech Sciences Foundation (GA\v CR), Grant Agreement
		24--11034S. The Institute of Mathematics of the Academy of Sciences of
		the Czech Republic is supported by RVO:67985840. } }

\date{}

\maketitle

{\small

\centerline{BCAM -- Basque Center for Applied Mathematics}
\centerline{Alameda de Mazarredo 14, E-48009 Bilbao, Basque Country, Spain}
\medbreak
\centerline{Ikerbasque -- Basque Foundation for Science}
\centerline{Plaza Euskadi 5, E-48009 Bilbao, Basque Country, Spain}
\medbreak
\centerline{Universit\'e Claude Bernard Lyon 1, ICJ UMR5208,}
\centerline{F-69622 Villeurbanne, France} %Institut Camille Jordan UMR CNRS 5208, Universit\'e Claude Bernard Lyon 1;}
%\centerline{43, Boulevard du 11 novembre 1918, F-69622 Villeurbanne, France}
\medbreak
\centerline{Email address: \texttt{ffanelli@bcamath.org}}

\bigbreak
\bigbreak
\centerline{Institute of Mathematics of the Academy of Sciences of the Czech Republic}
\centerline{\v Zitn\'a 25, CZ-115 67 Praha 1, Czech Republic}
\medbreak
\centerline{Email address: \texttt{feireisl@math.cas.cz}}

}

\medbreak

\begin{abstract}
	
	In contrast with a large variety of conventional models of thermally driven fluids, we show that the standard Oberbeck--Boussinesq approximation \emph{cannot} be obtained as a singular limit of the Navier--Stokes--Fourier system in the rotational coordinate system, with the buoyancy force proportional to the sum of the gravitational and centrifugal forces multiplied by the temperature variation.

\end{abstract}

%\bigskip

{\small


\medbreak

\tableofcontents

}

\section{Introduction}
\label{i}

The Oberbeck--Boussinesq (OB) approximation is widely used as a simple model of a fluid driven by thermally induced convection. There are numerous  
studies concerning mostly formal derivation of the model from the primitive Navier--Stokes--Fourier (NSF) system of a compressible, viscous, and heat conductive fluid, see e.g. the survey of Zeytounian \cite{ZEY1} and the references therein.
The formal argument giving rise to the OB approximation consists in
keeping the fluid density constant except in the gravitational force, where 
it is replaced by the (negative) temperature deviation. This process produces 
a non--potential driving force acting in the direction opposite to gravitation
in the momentum equation. 
The resulting system of equations is then supplemented by the conventional heat equation for the temperature.

A rigorous derivation of the OB system as an incompressible limit of the NSF system was performed in \cite[Chapter 5]{FeNo6A}. As illustrated in a series 
of recent studies \cite{BaBeFeOsTi, BelFeiOsch, Fa-Fe_2024}, the specific form of the limit system may depend not only on the underlying field equations, but also on the boundary conditions. In the context of the Rayleigh--B\' enard convection problem, the Dirichlet boundary conditions imposed on the temperature in the primitive system result in a rather unexpected \emph{non--local} boundary term in the limit system, see \cite{BelFeiOsch, Fa-Fe_2024}. Similarly, thermally driven compressible fluids under strong stratification produce a reduced Majda type system (Majda \cite{MAJ2})  rather than the conventionally used anelastic approximation, see \cite{BaBeFeOsTi}.

The goal of this work is to identify the singular limit of the NSF system in the incompressible, stratified regime in a rotating frame. In this context,  
the OB approximation is widely used, for instance, in models of geophysical flows, see Ecke and Shishkina \cite{EckeS} and the references therein;
see also Arslan \cite{Ars}, Kannan and Zhu \cite{K-Z}, and Welter \cite{Welter}, where (neglecting effects due to the centrifugal force)
the authors study bounds on the Nusselt number for rotating Rayleigh--B\'enard convection. The leading idea, proposed in the 
seminal work of Chandrasekhar \cite{CHAN}, consists in augmenting the buoyancy force 
in the standard OB approximation by a component resulting from the action of 
the centrifugal force in the rotating frame. As reported in \cite{EckeS}, see also Becker et al. \cite{BSCA}, the numerical computations seem to be in a good agreement with experiments.

In contrast with the common belief in the validity of the modified OB approximation for rotating fluids, we show that the limit system containing an active contribution from the temperature augmented centrifugal force is actually very \emph{different}. In particular, the fluid motion is purely horizontal, entirely 
independent of the vertical coordinate. Heuristically, we argue as follows:
\begin{itemize}
\item The origin of the Rayleigh--B\' enard convection flow is due to the \emph{compressibility} of the fluid. In the incompressible limit, where the pressure 
becomes constant, the density variation may be replaced by temperature variation with the opposite sign. This is the celebrated Boussinesq relation.  

\item For the limit system to feel this change, the Mach and Froude characteristic numbers must be properly scaled, cf. \cite{FeNo6A}. 

\item If the limit system is influenced by the centrifugal force, the same scaling must be applied. However, the scaled centrifugal force imposes imperatively 
smallnes of the Rossby number in the Coriolis force. 

\item In the regime of small Rossby number, the Taylor--Proudman theorem applies, enforcing the motion of the fluid to become purely horizontal.
This is in contrast with the conventional OB approximation, where a strong vertical movement results from the competition 
between the buoyancy and gravitational forces.

\end{itemize}
In order to state rigorously our results, we start with the precise formulation of the problem.

%Indeed the Boussinesq relation reversing the density and temperature variations in the incompressible limit requires the same singular scaling of the gravitational and centrifugal forces.
%\fra{\fbox{Remark:} i do not know if it is really necessary, generally speaking, to have the same singular scaling (i guess this depends
%on the type of regime one wants to consider?). However, what is sure is that this is precisely the scaling considered in \cite{EckeS}.}
%This in turn produces singularity in the 
%Coriolis force eliminating the vertical movement of the fluid in the asymptotic limit. The resulting fluid motion is necessarily planar in contrast with the conventional OB approximation.   

\subsection{Problem formulation} \label{ss:problem-formul}
We formulate the problem in the geometric framework considered in the review paper by Ecke and Shishkina \cite{EckeS}.

\subsubsection{Physical domain}

We suppose that the physical domain $\Omega \subset R^3$ occupied by the fluid is a cylinder 
\begin{equation} \label{s1}
\Omega = B(r) \times (0,1),\qquad  B(r) = \left\{ \vc x_h = (x_1, x_2) \ \Big| \ |\vc x_h | < r \right\}.  
\end{equation}
Here and hereafter, we write $\vc x\in\Omega$ in the form $\vc x = (\vc x_h, x_3)$ to stress the anisotropy between ``horizontal'' and ``vertical'' variables,
$\vc x_h=(x_1,x_2)$ and $x_3$, respectively. In addition, 
we assume the fluid domain is rotating around its vertical axis. Accordingly, it is convenient 
to write the equations of motion in the rotating coordinate frame. 
As a matter of fact, more general \emph{simply} connected rotating 2d domains could be considered.

\subsubsection{Primitive Navier--Stokes--Fourier system} \label{sss:equations}

Let $\vr = \vr(t,x)$, $\vt = \vt(t,x)$, and $\vu = \vu(t,x)$ denote the fluid mass density, the (absolute) temperature, and the velocity field, 
respectively.
The classical principles of conservation of mass, linear momentum, and energy written in the rotating frame give rise to the following 
(scaled) 
Navier--Stokes--Fourier (NSF) system: 
 
\begin{align} 
	\partial_t \vr + \Div (\vr \vu) &= 0, \label{eq:mass}\\
	\partial_t (\vr \vu) + \Div (\vr \vu \otimes \vu) + \frac{1}{\ep^2} \Grad p(\vr, \vt) + \frac{2}{\sqrt{\veps}} \vc{e}_3\times \vr \vu
	&= \Div \mathbb{S}(\vt, \Ds \vu) + \frac{1}{\ep} \vr \Grad G  + \frac{1}{2\ep} \vr \Grad |\vc{x}_h |^2 , \label{eq:mom} \\ 
	\partial_t (\vr e(\vr, \vt)) + \Div (\vr e(\vr, \vt) \vu) + \Div  \vc{q} (\vt, \Grad \vt) &= 
	\ep^2 \mathbb{S} (\vt, \Ds \vu) : \Ds \vu - p(\vr, \vt) \Div \vu .
	\label{eq:int-en}	
\end{align}

\noindent
Here, the functions $p(\vr,\vt)$ and $e(\vr,\vt)$ represent the pressure and the internal energy, respectively.
Their structural properties enforced by appropriate equations of state
will be specified in Section \ref{ss:constitutive}.
The function $G$
represents the gravitational potential acting in the vertical direction,
\begin{equation} \label{s5}
	G = - gx_3 .
\end{equation}
However, general gravitational fields given as 
\begin{equation*} 
G = \vc{g} \cdot x ,\ \vc{g} \in R^3,
\end{equation*}	
can be also considered.

The effect of rotation is represented by the Coriolis force $\vc{e}_3\times \vr\vu$, where $\vc{e}_3=(0,0,1)$, with the associated centrifugal 
force  $\vr \Grad|\vc{x}_h|^2$.

The viscous stress tensor is given by Newton's rheological law 
\begin{equation} \label{s6}
	\mathbb{S}(\vt, \Ds \vu) = 2\mu(\vt) \left( \Ds\vu - \frac{1}{3} \Div \vu \mathbb{I} \right) + \eta(\vt) \Div \vu \mathbb{I}, 
\end{equation}
where $\Ds\vu = \frac{1}{2}(\Grad \vu + \Grad^t \vu)$ is the symmetric part of the velocity gradient $\Grad \vu$. 
The heat flux is given by Fourier's law
\begin{equation} \label{s7} 
	\vc{q}(\vt, \Grad \vt) = - \kappa (\vt) \Grad \vt.
\end{equation}
The boundary and the initial data for system \eqref{eq:mass}--\eqref{eq:int-en} will be specified in Section \ref{s:primitive}.

System \eqref{eq:mass}--\eqref{eq:int-en} contains a (small) parameter $\ep > 0$ representing 
various kinds of scaling:
\begin{itemize}
\item 	The fluid is nearly incompressible, with the Mach number proportional to $\ep$.
\item   As shown in \cite[Chapter 5]{FeNo6A}, the appropriate scaling 
that gives rise to the Boussinesq relation between the density and temperature deviations in the buoyancy force is the Froude number of order $\sqrt{\ep}$.
This corresponds to the scaling of the gravitational force
\[
\frac{1}{\ep} \vr \Grad G.
\]
\item Anticipating a balance between the gravitation and the effect of the centrifugal force advocated in \cite{EckeS}, we suppose
\begin{equation} \label{s8}
	\mbox{centrifugal force}\ \approx \ \frac{1}{2\ep} \vr \Grad |\vc{x}_h |^2.
\end{equation}
\item At the same time, the scaling \eqref{s8} yields \emph{imperatively} the Rossby number in the Coriolis force proportional to $\sqrt{\ep}$,
\[
\frac{2}{\sqrt{\ep}} \vc{e}_3 \times \vr \vu. 
\]
\end{itemize}	

\noindent
Note carefully that the scaling of the Coriolis force is \emph{enforced} 
by $\eqref{s8}$. As we show below, it is essentially this fact that eliminates the standard OB approximation as 
a possible singular limit.

\subsection{Singular limit} \label{ss:invalid}

We are ready to discuss our claim that the  singular limit 
of NS system \eqref{eq:mass}--\eqref{eq:int-en} in the regime $\ep\to 0$ is rather different from the expected
(and commonly used) OB approximation.

At the leading order, the dynamics is driven by the pressure, the gravitation and the centrifugal force terms, all proportional to $\frac{1}{\ep}$.
As already observed in many standard situations, see e.g. \cite[Chapter 5]{FeNo6A} or \cite{DS-F-S-WK} for the case of rotating fluids,
this balance gives rise to the Boussinesq relation between the 
mass density and the temperature deviations, here augmented by an additional term owing to the scaling
of the centrifugal force.

At the first glance, the influence of the Coriolis force scaled as $\frac{1}{\sqrt{\ep}}$ seems negligible. Still
its effect can be captured by projecting the momentum equation onto the space of solenoidal vector fields, yielding non-trivial constraints on the target dynamics,
see \cite{DS-F-S-WK, FeGaGVNo, FeNo8} among others. In particular,
the fast rotation produces a vertical rigidity, known in the geophysical realms as \emph{Taylor--Proudman theorem}. 
Indeed, the Coriolis force eliminates entirely the vertical motion from
the asymptotic dynamics.

In view of the above arguments, the fluid motion in the asymptotic limit $\ep \to 0$ becomes
necessarily \emph{planar}, in sharp contrast with the conventional OB approximation. This is even more surprising
in the context of the Rayleigh-B\'enard problem, where the vertical motion is enhanced not only by  the gravitation but also by the 
strong buoyancy force caused by a non-zero background temperature gradient acting from the bottom to the top boundary of the domain.

\medbreak
The goal of the paper is to provide a rigorous justification of the above heuristic arguments. We start by introducing the main 
hypotheses and basic properties of the NSF system in Section \ref{s:primitive}. The main results are then stated in Section \ref{s:result}, see Theorem \ref{TM1}. 
In Section \ref{u}, we derive the necessary uniform bounds on the family of scaled solutions, independent of the parameter $\ep \to 0$.
In Section \ref{A}, we characterise the asymptotic dynamics in the limit $\ep\to0$, thus completing the proof of the main results.

\section{Mathematics of the Navier--Stokes--Fourier system} 
\label{s:primitive}

In this section, we introduce the basic hypotheses imposed on constitutive relations and recall some well known facts 
concerning the primitive NSF system.

\subsection{Constitutive relations} \label{ss:constitutive}

To close system \eqref{eq:mass}--\eqref{eq:int-en}, we have to specify the constitutive relations, namely the equations of state (EOS) and 
the form of transport coefficients.
These are similar to \cite[Chapters 1,2]{FeNo6A} (see also \cite[Chapter 1]{FeiNovOpen}), and they are motivated by the available \emph{existence} theory. 

The pressure EOS reads
\[
p(\vr, \vt) = p_{\rm m} (\vr, \vt) + p_{\rm rad}(\vt), 
\]
where $p_{\rm m}$ is the pressure of a general \emph{monoatomic} gas related to the associated internal energy as 
\begin{equation} \label{con1}
	p_{\rm m} (\vr, \vt) = \frac{2}{3} \vr e_{\rm m}(\vr, \vt).
\end{equation}
The symbol $p_{\rm rad}$ is the radiation pressure, which takes the form
\[
p_{\rm rad}(\vt) = \frac{a}{3} \vt^4,\qquad a > 0.
\]
The radiation pressure plays a crucial role in the existence theory  
developed in \cite{DF1, FeNo6A} eliminating possible uncontrolled temperature oscillations 
in the (hypothetical) vacuum zones.
The pressure $p_m$ can be more general in the sense specified in \cite[Chapter 1, Section 1.4]{FeNo6A}.

The pressure and the internal energy are interrelated through the Gibbs law
\begin{equation} \label{Gibbs}
\vt Ds = De + p D \left( \frac{1}{\vr} \right),
\end{equation}
where $s$ is the (specific) entropy. In addition, we impose the hypothesis of thermodynamic stability
\begin{equation} \label{ThSt}
\frac{\partial p(\vr, \vt)}{\partial \vr} > 0,\ \frac{\partial e(\vr, \vt)}{\partial \vt} > 0. 
\end{equation} 

%The specific form of \eqref{con1} is another issue to be discussed. In the context of gases, 
%the natural candidate is provided by the standard Boyle--Mariotte law $p = \vr \vt$, $e = \frac{3}{2} \vt$. Unfortunately, this \textsl{ansatz} 
%fails to provide even the expected energy estimates as long as the boundary temperature is prescribed. To circumvent 
%this difficulty, more physics must be taken into account, cf. \cite[Chapter 1]{FeNo6A}, as specified here below.

Our main hypotheses concerning the EOS are formulated below:
\begin{itemize}
	
	\item \tbf{Gibbs' law} together with \eqref{con1} yield 
	\[
	p_{\rm m} (\vr, \vt) = \vt^{\frac{5}{2}} P \left( \frac{\vr}{\vt^{\frac{3}{2}}  } \right),
	\]
	for a certain $P \in C^1[0,\infty)$.
	Consequently, 
	\begin{equation} \label{w9}
		p(\vr, \vt) = \vt^{\frac{5}{2}} P \left( \frac{\vr}{\vt^{\frac{3}{2}}  } \right) + \frac{a}{3} \vt^4,\quad
		e(\vr, \vt) = \frac{3}{2} \frac{\vt^{\frac{5}{2}} }{\vr} P \left( \frac{\vr}{\vt^{\frac{3}{2}}  } \right) + \frac{a}{\vr} \vt^4, \qquad a > 0.
	\end{equation}
	
	\item The \tbf{hypothesis of thermodynamic stability} %\eqref{HTS}
	expressed in terms of  $P$ gives rise to
	\begin{equation} \label{w10}
	P \in C^1[0, \infty),\ P(0) = 0, \ P'(Z) > 0 \ \mbox{for}\ Z \geq 0,\qquad 0<  \frac{ \frac{5}{3} P(Z) - P'(Z) Z }{Z} \leq c \ \mbox{ for }\ Z > 0,
	\end{equation} 	
where the upper bound in the last condition means boundedness of the specific heat at constant volume.	
		\item 
	The associated specific {\bf entropy} takes the form 
	\begin{equation} \label{w12}
		s(\vr, \vt) = s_{\rm m}(\vr, \vt) + s_{\rm rad}(\vr, \vt),\qquad s_{\rm m} (\vr, \vt) = \mathcal{S} \left( \frac{\vr}{\vt^{\frac{3}{2}} } \right),\qquad
		s_{\rm rad}(\vr, \vt) = \frac{4a}{3} \frac{\vt^3}{\vr}, 
	\end{equation}
	where 
	\begin{equation} \label{w13}
		\mathcal{S}'(Z) = -\frac{3}{2} \frac{ \frac{5}{3} P(Z) - P'(Z) Z }{Z^2} < 0.
	\end{equation}

\item {\bf Third law of thermodynamics}: motivated by \cite{FeiLuSun}, we impose the {Third law of thermodynamics}, cf.~Belgiorno \cite{BEL1, BEL2}. 
Specifically, we require  the entropy to vanish 
when the absolute temperature approaches zero, 
\begin{equation} \label{w14}
	\lim_{Z \to \infty} \mathcal{S}(Z) = 0.
\end{equation}
In addition, we suppose 
\begin{equation} \label{w14a}
P \in C^1[0,\infty) \ \mbox{ is such that } \ \liminf_{Z\to\infty}\frac{P(Z)}{Z}>0, 
\end{equation}
see Section 2.1.1 of \cite{FeiLuSun} for details.

\end{itemize}

It is interesting to note that all the above restrictions except \eqref{w14} 
imposed on $p_m$ are also satisfied by the conventional Boyle--Mariotte law corresponding to $P(Z) = Z$. Moreover, as shown in \cite[Section 2.2.1]{FeiLuSun}, 
the hypotheses \eqref{w14}, \eqref{w14a} yield coercivity of the pressure law, specifically, the function $Z \mapsto P(Z)/ Z^{\frac{5}{3}}$ is decreasing, and  
\begin{equation} \label{w11}
	\lim_{Z \to \infty} \frac{ P(Z) }{Z^{\frac{5}{3}}} = p_\infty > 0.
\end{equation}

As for the transport coefficients, we suppose they are continuously differentiable functions of the 
temperature satisfying
\begin{align} 
	0 < \underline{\mu}(1 + \vt) &\leq \mu(\vt),\qquad |\mu'(\vt)| \leq \Ov{\mu}, \br 
	0 &\leq \underline{\eta} (1 + \vt) \leq \eta (\vt) \leq \Ov{\eta}(1 + \vt), \br
	0 < \underline{\kappa} (1 + \vt^\beta) &\leq \kappa (\vt) \leq \Ov{\kappa}(1 + \vt^\beta), 
	\quad \mbox{ where }\ \beta > 6. \label{w16}
\end{align}
Similarly to the hypotheses imposed on EOS, the 
restriction $\beta > 6$ is dictated by the available existence theory, cf. \cite{FeiNovOpen}.

As a consequence of the hypotheses \eqref{w9} -- \eqref{w14a}, we get the following bounds:
\begin{align} 
	\vr^{\frac{5}{3}} + \vt^4 \lesssim \vr e(\vr, \vt) &\lesssim 	1+ \vr^{\frac{5}{3}} + \vt^4, \label{L5b} \\
	s_{\rm m}(\vr, \vt) &\lesssim \left( 1 + |\log(\vr)| + [\log(\vt)]^+ \right), \label{L5a}
\end{align} 
see \cite[Chapter 3, Section 3.2]{FeNo6A}.  

\subsection{Boundary and initial conditions}

We impose the conventional no-slip boundary conditions on the lateral 
boundary of the cylinder,
\begin{equation} \label{B1}
\vu|_{ \partial B(r) \times (0,1) } = 0, 
\end{equation}
supplemented with the complete--slip at the top and bottom parts, 
\begin{equation} \label{B2}
\vu \cdot \vc{n}|_{ B(r) \times \{ x_3 = 0, 1 \} } = 0,\ 
[\mathbb{S}(\vt, \Ds \vu) \cdot \vc{n} ] \times \vc{n} |_{ B(r) \times \{ x_3 = 0, 1 \} } = 0.
\end{equation}

\begin{Remark} \label{r:BL}
In the context of fast rotating fluids, it is well-known that no-slip boundary conditions give rise to boundary layer phenomena, the so--called Ekman boundary layers
(in proximity of horizontal boundaries) and Munk boundary layers (near vertical walls). None of these effects will actually appear in our study.
We postpone comments about this issue at the end of the paper, see Section \ref{s:comments}.
\end{Remark}

In accordance with the given scaling, we may consider the 
Dirichlet boundary conditions for the temperature in the form
\begin{align} 
	\vt_{\ep,B} &= \Ov{\vt} + \ep \mathfrak{T}_{\rm bot}  \quad \mbox{ if }\ x_3 = 0	, \ |\vc{x}_h| \leq  r, \br 
	\vt_{\ep,B} &= \Ov{\vt} + \ep \mathfrak{T}_{\rm top}  \quad \mbox{ if }\ x_3 = 1,\ |\vc{x}_h| \leq r, \br 
	\vt_{\ep,B} &= \Ov{\vt} + \ep x_3 \mathfrak{T}_{\rm top} + \ep (1 - x_3) \mathfrak{T}_{\rm bot}\quad \mbox{ if }\ |\vc{x}_h| = r ,
	\nonumber
\end{align}
where $\Ov{\vt} > 0$ is a constant,  and $\mfr T_\bot$, $\mfr T_\top$ are smooth functions defined on $R^2$. We can actually handle a more general setting 
\begin{equation} \label{eq:t_boundary}
\vt|_{\partial \Omega} = \vt_{\ep, B},\ 
\vt_{\ep, B} = \Ov{\vt} + \ep \vtB, \ \Ov{\vt} > 0 \ \mbox{constant},
\end{equation}	
where $\vtB$ is a restriction of a smooth function defined on $R^3$.

Finally, anticipating again the chosen scaling, we consider the initial conditions in the form 
\begin{align} 
\vr(0, \cdot) &= \vr_{0, \ep} = \Ov{\vr} + \ep \mathcal{R}_{0, \ep}, \br	
\vt(0, \cdot) &= \vt_{0,\ep} = \Ov{\vt} + \ep \mathfrak{T}_{0, \ep},
\label{B3}	
\end{align}
where $\Ov{\vt}$ is the constant introduced in the boundary condition \eqref{eq:t_boundary}, and $\Ov{\vr} > 0$ is another constant chosen so that 
\[
\avintO{\vr_{0,\ep}} = \Ov{\vr} ,
\ \mbox{ meaning }\ \intO{ \mathcal{R}_{0, \ep} } = 0.
\]
Here and hereafter, the symbol $\avintO{f} \equiv \frac{1}{|\Omega|} \intO{f}$ stands for the integral average over $\Omega$.
In accordance with the boundary conditions \eqref{B1}, \eqref{B2}, the boundary is impermeable, in particular, the total mass 
\[
\intO{ \vr(t,\cdot) } 
\]
is a constant of motion.

\subsection{Weak formulation of the NSF system} \label{ss:weak}

The weak formulation of the primitive NSF system follows the leading idea proposed in \cite[Chapter 3]{FeNo6A}, namely replacing 
the internal energy balance \eqref{eq:int-en} by the entropy inequality supplemented by the total energy balance. 
Later \cite{ChauFei, FeiNovOpen}, the total energy was replaced by the ballistic energy to accommodate the Dirichlet boundary conditions for the temperature.

\begin{Definition}[\tbf{Weak solution to NSF system}] \label{DL1}
	We say that a trio $(\vre, \vte, \vue)$ is a \emph{weak solution} of the scaled NSF system \eqref{eq:mass}--\eqref{eq:int-en},
	supplemented with 
	%the constitutive relations of Subsection \ref{ss:constitutive},
	the boundary conditions \eqref{B1}, \eqref{B2}, \eqref{eq:t_boundary}, and
	the initial data
	\[
\vre(0, \cdot) = \vr_{0,\ep},\ \vte(0, \ep) = \vt_{0,\ep},\ 
	\vue(0, \cdot) = \vu_{0,\ep},
	\]
	if the following holds true:
	\begin{itemize}
		
		\item The solution belongs to the {\bf regularity class} 
		\begin{align}
			\vre &\in L^\infty(0,T; L^{\frac{5}{3}}(\Omega)),\ \vre \geq 0 
			\ \mbox{a.a.~in}\ (0,T) \times \Omega, \br
			\vue &\in L^2(0,T; W^{1,2} (\Omega; R^3)), \ 
			\vue |_{\partial B(r) \times [0,1]} = 0,\  
			\vue \cdot \vc{n}|_{x_3 = 0,1} = 0, \br 
			\vte^{\beta/2} ,\ \log(\vte) &\in L^2(0,T; W^{1,2}(\Omega)) \ \mbox{for some}\ \beta \geq 2,\ 
			\vt > 0 \ \mbox{a.a.~in}\ (0,T) \times \Omega, \br
			(\vte - \vt_{\ep,B}) &\in L^2(0,T; W^{1,2}_0 (\Omega)).
			\label{Lw6}
		\end{align}
		
		\item The {\bf equation of continuity} \eqref{eq:mass} is satisfied in the sense of distributions including the impermeability boundary conditions, specifically
		\begin{align} 
			\int_0^T \intO{ \Big[ \vre \partial_t \varphi + \vre \vue \cdot \Grad \varphi \Big] } \dt &=  - 
			\intO{ \vr_{0,\ep} \varphi(0, \cdot) }
			\label{Lw4}
		\end{align}
		for any $\varphi \in C^1_c([0,T) \times \Ov{\Omega} )$.
		\item The {\bf momentum equation} \eqref{eq:mom} is satisfied in the sense of distributions, specifically 
		\begin{align}
			\int_0^T &\intO{ \left[ \vre \vue \cdot \partial_t \bfphi + \vre \vue \otimes \vue : \Grad \bfphi - \frac{2}{\sqrt{\ep}} (\mbf{e}_3 \times \vre \vue) \cdot 
				\bfphi + 
				\frac{1}{\ep^2} p(\vre, \vte) \Div \bfphi \right] } \dt \br &= \int_0^T \intO{ \left[ \mathbb{S}(\vte, \Ds \vue) : \Grad \bfphi - \frac{1}{\ep} \vre \Grad G \cdot \bfphi -  \frac{1}{2\ep} \vre \Grad |\vc{x}_h |^2 \cdot \bfphi  \right] } \dt \br &- 
			\intO{ \vr_{0,\ep} \vu_{0,\ep} \cdot \bfphi (0, \cdot) }
			\label{Lw5}
		\end{align}	
		for any $\bfphi \in C^1_c([0, T) \times \Ov{\Omega}; R^3)$ such that 
		$\bfphi \cdot \vc{n}|_{\partial \Omega} = 0$.
		
		\item  The {\bf entropy balance} is satisfied as inequality 
		\begin{align}
			- \int_0^T &\intO{ \left[ \vre s(\vre, \vte) \partial_t \varphi + \vre s (\vre ,\vte) \vue \cdot \Grad \varphi + \frac{\vc{q} (\vte, \Grad \vte )}{\vte} \cdot 
				\Grad \varphi \right] } \dt \br &\geq \int_0^T \int_{\Omega}{\frac{\varphi}{\vte} \left( \ep^2 \mathbb{S}(\vte, \Ds \vue) : \Ds \vue - \frac{\vc{q}(\vte ,\Grad \vte)\cdot \Grad \vte}{\vte}\right)} \dx \dt  \br &+ \intO{ \vr_{0, \ep} s(\vr_{0, \ep}, \vt_{0,\ep}) 
				\varphi (0, \cdot) } 
			\label{Lw7} 
		\end{align}
		for any $\varphi \in C^1_c([0, T) \times \Omega)$, $\varphi \geq 0$. 
		
		\item  The {\bf ballistic energy balance}
		\begin{align}  
			- &\int_0^T \partial_t \psi	\intO{ \left[ \ep^2 \frac{1}{2} \vre |\vue|^2 + \vre e(\vre, \vte) - \wtilde \vt_{\ep,B} \vre s(\vre, \vte) \right] } \dt  \br &+ \int_0^T \int_{{\Omega}} \psi \frac{\wtilde \vt_{\ep,B} }{\vte} \left( \ep^2 \mathbb{S}(\vte, \Ds \vue) : \Ds \vue - \frac{\vc{q}(\vte ,\Grad \vte)\cdot \Grad \vte}{\vte}\right)   \dx \dt 
			\br
			&\leq 
			\int_0^T \psi \intO{ \left[ \ep \vre \vue \cdot \Grad G +  \frac{\ep}{2} \vre \vue \cdot \Grad |\vc{x}_h |^2 \right] } \dt \br  &- 
			\int_0^T \psi \intO{ \left[ 
				\vre s(\vre, \vte) \partial_t \wtilde \vt_{\ep,B} + \vre s(\vre, \vte) \vue \cdot \Grad \wtilde \vt_{\ep,B}  + \frac{\vc{q}(\vte, \Grad \vte)}{\vte}\cdot \Grad \wtilde \vt_{\ep,B} \right] } \dt \br 
			&+ \psi(0) \intO{ \left[ \frac{1}{2} \ep^2 \vr_{0,\ep} |\vu_{0,\ep}|^2 + \vr_{0,\ep} e(\vr_{0,\ep}, \vt_{0, \ep}) - \wtilde \vt_{\ep,B} (0, \cdot) \vr_{0,\ep} s(\vr_{0,\ep}, \vt_{0,\ep}) \right] }
			\label{Lw8}
		\end{align}
		holds true for any $\psi \in C^1_c ([0, T))$, $\psi \geq 0$, and \emph{any} continuously differentiable extension $\wtilde \vt_{\ep,B}$ of the boundary datum, 
		\[
		\wtilde \vt_{\ep,B} > 0 \ \mbox{in}\ [0,T] \times \Ov{\Omega},\ \wtilde \vt_{\ep,B}|_{\partial \Omega} = \vt_{\ep, B}.
		\]
	\end{itemize}

\end{Definition}

The \emph{existence} of global in time weak solutions in the sense of Definition \ref{DL1} for the no-slip boundary conditions 
for the velocity was shown in \cite{ChauFei}, \cite[Chapter 12]{FeiNovOpen} on condition that $\Omega$ is a smooth (at least 
$C^2$) domain. 
The proof can be easily modified to accommodate the present mixed boundary conditions for the velocity. In general, the lack of smoothness of the domain can be an issue, however, the present cylindrical shape  
can be accommodated in the existence proof as long we can construct
a harmonic extension of the boundary data of class $C^1$ (cf. \cite[Chapter 12, Section 12.4.1]{FeiNovOpen}). Specifically, 
to obtain the necessary uniform bounds, we need a function $\widetilde{\vt}_B$,
\[
\Del \widetilde{\vt}_B = 0,\ 
\widetilde{\vt}_B|_{\partial \Omega} = \vtB,
\]
with a bounded gradient. The problem can be rewritten in the form 
\[
\widetilde{\vt}_B = \xi + \vtB,\ 
\Del \xi = - \Del \vtB,\ \xi|_{\partial \Omega} = 0.
\] 
As $\vtB$ is of class $C^2$ we can extend $\Del \vtB$ as a 2-periodic 
odd function in the $x_3$ variable preserving the property\footnote{Here, the notation $[-1,1]|_{\{ - 1, 1\} }$ stands for the one-dimensional
torus constructed over the interval $[-1,1]$, after identification of the points $-1$ and $1$.} 
$\Del \vtB \in L^p(B(r) \times [-1,1]|_{\{ - 1, 1\} })$ 
for any finite $p$. Using the standard elliptic estimates on the 
periodic domain $B(r) \times [-1,1]|_{\{ - 1, 1\} }$ we conclude 
$\xi \in W^{2,p}(\Omega)$ for any $1 \leq p < \infty$. In particular, the extension $\widetilde{\vt}_B$ admits a bounded 
gradient as long as $p > 3$.

\section{Main result} \label{s:result}

Before stating our main result, let us introduce several material parameters:
\begin{itemize}
 \item the thermal expansion coefficient
\[
\alpha(\Ov{\vr}, \Ov{\vt} ) \equiv \frac{1}{\Ov{\vr}}  \frac{\partial p(\Ov{\vr}, \Ov{\vt} ) }{\partial \vt} \left( \frac{\partial p(\Ov{\vr}, \Ov{\vt} ) }{\partial \vr} \right)^{-1};
\]
\item the specific heat at constant pressure and constant volume
\[
 c_p (\Ov{\vr}, \Ov{\vt} ) \equiv \frac{\partial e(\Ov{\vr}, \Ov{\vt} ) }{\partial \vt}	+ \Ov{\vr}^{-1} \Ov{\vt}
	\alpha(\Ov{\vr}, \Ov{\vt} ) \frac{\partial p(\Ov{\vr}, \Ov{\vt} ) }{\partial \vt},\ c_v (\Ov{\vr}, \Ov{\vt} ) \equiv \frac{\partial e(\Ov{\vr}, \Ov{\vt} ) }{\partial \vt};
\]
\item the coefficient
\[
\lambda(\Ov{\vr}, \Ov{\vt}) \equiv \frac{\Ov{\vt} \alpha (\Ov{\vr}, \Ov{\vt} ) }
{\Ov{\vr} c_p(\Ov{\vr}, \Ov{\vt}) } \frac{\partial p(\Ov{\vr}, \Ov{\vt} )}{\partial \vt} = 1 - \frac{c_v(\Ov{\vr}, \Ov{\vt}) }{c_p(\Ov{\vr}, \Ov{\vt} ) } \ \in (0,1).
\] 
\end{itemize}

We are ready to formulate our main result.

\begin{Theorem}[{\bf Singular limit $\ep \to 0$}] \label{TM1}
	
Let $(\vre, \vte, \vue)_{\ep > 0}$ be a family of weak solutions 
of the scaled NSF system with the boundary conditions 
\eqref{B1}--\eqref{eq:t_boundary}, emanating from the initial data 
\[
\vre(0, \cdot) = \Ov{\vr} + \ep \mathcal{R}_{0,\ep} ,\ \vue(0, \cdot) =  \vu_{0, \ep}, \ 
\vte(0, \cdot) = \Ov{\vt} + \ep \mathfrak{T}_{0, \ep}. 
\]
In addition, suppose 
\begin{align} 
\| \mathcal{R}_{0, \ep} \|_{L^\infty(\Omega)} \aleq 1,\ 
\intO{ \mathcal{R}_{0,\ep} } = 0,\ \mathcal{R}_{0,\ep} &\to \mathcal{R}_0 \ \mbox{in}\ L^1(\Omega); \br 
\| {\mathfrak{T}}_{0, \ep} \|_{L^\infty(\Omega)} \aleq 1,\ 
 \mathfrak{T}_{0,\ep} &\to \mathfrak{T}_0 \ \mbox{in}\ L^1(\Omega), \br
\| {\vc{u}}_{0, \ep} \|_{L^\infty(\Omega;R^3)} \aleq 1,\ 
\vc{u}_{0,\ep} &\to \vc{u}_0 \ \mbox{in}\ L^1(\Omega;R^3),\ \mbox{where}\ 
\vu_0 = [\vu_{0,h}, 0],\ \vu_{0,h} = \vu_{0,h}(\vc{x}_h).
\label{LLw1}
\end{align}	 
Suppose also that
\begin{align} 
\mathfrak{T}_0 &\in W^{2,p}(\Omega) \ \mbox{for all}\ 1 \leq p < \infty, \ \mathfrak{T}_0|_{\partial \Omega} = \vtB, \br 
\vu_{0,h} &\in W^{2,p}(B(r); R^2) \ \mbox{for all}\ 1 \leq p < \infty,  
\vu_{0,h}|_{\partial B(r)} = 0,\ \Divh \vu_{0,h} = 0, 	
\label{LLw2}
\end{align}
and	that
\begin{equation} \label{LLw3}
\frac{\partial p(\Ov{\vr}, \Ov{\vt})}{\partial \vr } \Grad \mathcal{R}_0 + 
\frac{\partial p(\Ov{\vr}, \Ov{\vt})}{\partial \vt } \Grad \mathfrak{T}_0 = \Ov{\vr} \Grad \left( G +  |x_h|^2 \right).	
\end{equation}	
	
Then
\begin{align} 
\frac{\vre - \Ov{\vr}}{\ep} &\to \mathcal{R} \ \mbox{in}\ L^\infty(0,T; L^1(\Omega)), \ \mbox{ with }\ \avintO{\mathcal{R}(t,\cdot)} = 0, \br	
\frac{\vte - \Ov{\vt}}{\ep} - 
\lambda(\Ov{\vr}, \Ov{\vt}) \avintO{\frac{\vte - \Ov{\vt}}{\ep}}   &\to \Theta \ \mbox{in}\ L^\infty(0,T; L^1(\Omega)) \cap 
L^2(0,T; W^{1,2}(\Omega)), \br
\vue &\to \vc{U} = (\vuh,0) \ \mbox{in}\ L^2(0,T; W^{1,2}(\Omega; R^3)),\ \mbox{where}\ \vu_h = \vuh (t, \vc{x}_h), \br  
\sqrt{ \vre } \vue &\to \sqrt{\Ov{\vr}} \ \vU \ % (\vuh, 0) \ 
\mbox{in}\ L^\infty(0,T; L^2(\Omega; R^3)) 
\nonumber
\end{align}
as $\ep \to 0$, where $(\mathcal{R}, \Theta, \vuh)$ is the (unique) strong solution of the target system (TS):
\begin{align}
\Divh \vuh &= 0, \br	
\Ov{\vr} \Big[ \partial_t \vuh +  \Divh (\vuh \otimes \vuh) \Big] + \Gradh \Pi &= \mu(\Ov{\vt}) \Delta_h \vuh + \left< \mathcal{R} \right> 
\Gradh\Big( G + \frac{1}{2}|\vc x_h|^2 \Big) ,\ \left< \mathcal{R} \right> \equiv \int_0^1 
\mathcal{R}(t,\vc x_n, x_3) \ \D x_3, \br 
\mbox{in}\ (0,T) &\times B(r), \br 
\label{TS1}
\end{align}	
\begin{align}
\Ov{\vr} c_p(\Ov{\vr}, \Ov{\vt}) \Big[ \partial_t \Theta + \vuh \cdot 
\Gradh \Theta \Big] - \Ov{\vr} \Ov{\vt} \alpha (\Ov{\vr}, \Ov{\vt} ) 
\vuh \cdot \Gradh \Big( G + \frac{1}{2}|\vc x_h|^2 \Big) &= \kappa (\Ov{\vt}) \Del \Theta , \br 
\mbox{in}\ (0,T) &\times \Omega,
\label{TS2}	
\end{align} 
supplemented with the Boussinesq relation 
\begin{equation} \label{TS3}
\frac{\partial p(\Ov{\vr}, \Ov{\vt})}{\partial \vr} \Grad \mathcal{R} 
+	\frac{\partial p(\Ov{\vr}, \Ov{\vt})}{\partial \vt} \Grad \Theta = 
\Ov{\vr} \Grad \Big( G + \frac{1}{2}|\vc x_h|^2 \Big),
\end{equation}		
the boundary conditions 
\begin{equation} \label{TS4}
\vuh|_{\partial B(r)} = 0, \ \Theta|_{\partial \Omega} = \vtB - 
\frac{\lambda(\Ov{\vr}, \Ov{\vt})}{ 1 -\lambda(\Ov{\vr}, \Ov{\vt}) } 
\avintO{ \Theta },	
\end{equation}		
and the initial conditions
\begin{equation} \label{TS5}
\vuh(0, \cdot) = \vu_{0,h},\ \Theta(0, \cdot) = \mathfrak{T}_0 - 
\lambda(\Ov{\vr}, \Ov{\vt}) \avintO{ \mathfrak{T}_0 }.	
\end{equation}	
\end{Theorem}	

\begin{Remark} \label{TR1}
	
A direct manipulation reveals that
\begin{equation} \label{TS5a}
\frac{\lambda(\Ov{\vr}, \Ov{\vt})}{ 1 -\lambda(\Ov{\vr}, \Ov{\vt}) } =
\frac{c_p (\Ov{\vr}, \Ov{\vt} )}{c_v  (\Ov{\vr}, \Ov{\vt} ) } - 1 .
\end{equation} 		

\end{Remark}

As claimed in the introductory part, the limit fluid motion is purely horizontal, in contrast with the 
commonly accepted OB dynamics. The non--local boundary term in \eqref{TS4} is pertinent to the 
Dirichlet boundary conditions imposed on the temperature and has been identified in \cite{BelFeiOsch}. 
Note that the above scenario seems the \emph{only compatible} with the incompressible limit as long as 
the effect of the centrifugal force is anticipated.

The rest of the paper is devoted to the proof of Theorem \ref{TM1}. The reader will have noticed that the initial data of the NSF system are 
\emph{well prepared}. Similarly to \cite{BelFeiOsch}, the proof leans 
on the new concept of ballistic energy inequality introduced in \cite{ChauFei}, \cite{FeiNovOpen}.

We finish this section by stating the relevant global existence result 
for the target system.

\subsection{Solvability of the target system}

As the target momentum equation \eqref{TS1} reduces to a variation of the 2d Navier-Stokes system,
it is plausible to expect global existence of strong solutions to the target problem. This is indeed the case as shown in \cite[Proposition 4.1]{FeGwSG25}.

\begin{Proposition}[{\bf Existence for the target system}] \label{TP1}
Suppose the initial data $\vu_{0,h}$, $\mathfrak{T}_0$ belong
to the class specified in \eqref{LLw2}. 

Then the target system \eqref{TS1}--\eqref{TS5} admits a unique regular solution in the class
\begin{align}
\vuh &\in L^p(0,T; W^{2,p}(B(r); R^2)), \ \partial_t \vuh \in L^p(0,T; L^{p}(B(r); R^2)), \br
\Theta &\in L^p(0,T; W^{2,p}(\Omega)), \ \partial_t \Theta \in L^p(0,T; L^{p}(\Omega))	 
\label{TS6}	
\end{align}	
for all $1 \leq p < \infty$.
\end{Proposition}

\section{Uniform bounds}
\label{u}

Our first goal is to derive uniform bounds on the sequence of solutions of the scaled NSF system, namely bounds which are independent of the parameter $\ep \to 0$.

\subsection{Relative energy}	

Similarly to \cite{Fa-Fe_2024}, we consider the scaled energy functional
\begin{equation} \label{u1}
	E_\ep (\vr, \vt, \vu) = \frac{1}{2} \vr |\vu|^2 + \frac{1}{\ep^2} 
\vr e(\vr, \vt),
\end{equation}	
together with the associated relative energy
\begin{align} 
E_\ep &\left( \vr, \vt, \vu \Big| \tvr, \tvt, \tvu \right) = \frac{1}{2} \vr |\vu - \tvu |^2 \br 
&+\frac{1}{\ep^2} \left[ \vr e(\vr, \vt) - \tvt \Big( \vr s(\vr , \vt) - \tvr s(\tvr, \tvt) \Big) - \left( e(\tvr, \tvt) - \tvt s(\tvr, \tvt) + \frac{p(\tvr, \tvt)}{\tvr} \right)(\vr - \tvr) - \tvr e(\tvr, \tvt)    \right].	
\label{u2}
	\end{align}
	
As shown in \cite{ChauFei}, \cite[Chapter 12]{FeiNovOpen}, any weak solution $(\vre, \vte, \vue)$ of NSF system in the sense of
Definition \ref{DL1} satisfies the relative energy inequality, 
 \begin{align}
	&\left[ \intO{ E_\ep \left(\vre, \vte, \vue \Big| \tvr, \tvt, \tvu \right) } \right]_{t = 0}^{t = \tau} \br 
	&+ \int_0^\tau  \int_{{\Omega}} \frac{\tvt }{\vte} \left(  \mathbb{S}(\vte, \Ds \vue) : \Ds \vue - \frac{1}{\ep^2} \frac{\vc{q}(\vte ,\Grad \vte)\cdot \Grad \vte}{\vte}\right)   \dx \dt \br 
	&\leq - \frac{1}{\ep^2} \int_0^\tau \intO{ \left( \vre (s(\vre, \vte) - s(\tvr, \tvt)) \partial_t \tvt + \vre (s(\vre,\vte) - s(\tvr, \tvt)) \vue \cdot \Grad \tvt \right) } \dt \br 
	&+ \frac{1}{\ep^2} \int_0^\tau \intO{ 
		\frac{\kappa (\vte) \Grad \vte}{\vte}  \cdot \Grad \tvt } \dt \br 
	&- \int_0^\tau \intO{ \Big[ \vre (\vue - \tvu) \otimes (\vue - \tvu) + \frac{1}{\ep^2} p(\vre, \vte) \mathbb{I} - \mathbb{S}(\vte, \Ds \vue) \Big] : \Ds \tvu } \dt   \br
	&+ \int_0^\tau \intO{ \vre \left[ \frac{1}{\ep} \Grad G + \frac{1}{2\ep} \Grad |\vc{x}_h |^2 
		- \frac{2}{\sqrt{\ep}} (\vc{e}_3 \times \vue)  - \partial_t \tvu - (\tvu \cdot \Grad) \tvu  \right] \cdot (\vue - \tvu) } \dt  \br 
	&+ \frac{1}{\ep^2} \int_0^\tau \intO{ \left[ \left( 1 - \frac{\vre}{\tvr} \right) \partial_t p(\tvr, \tvt) - \frac{\vre}{\tvr} \vue \cdot \Grad p(\tvr, \tvt) \right] } \dt
	\label{u3}
\end{align}
for a.a. $\tau > 0$ and any trio of continuously differentiable functions $(\tvr, \tvt, \tvu)$ satisfying
\begin{equation} \label{u4}
	\tvr > 0,\quad \tvt > 0,\quad \tvt|_{\partial \Omega} = \vt_{\ep,B}, \quad \tvu|_{\partial B(r) \times (0,1)} = 0,\quad \tvu \cdot \vc{n}|_{x_3 = 0,1} = 0.
\end{equation}

It is convenient to use the notation introduced in 
\cite{FeNo6A} distinguishing the ``essential'' and ``residual'' range of the 
thermostatic variables $(\vr, \vt)$. Specifically, given a compact set 
\begin{equation*}
	K \subset \left\{ (\vr, \vt) \in \R^2 \ \Big| \ \vr > 0, \vt > 0 \right\}
\end{equation*}
and $\ep > 0$, we denote 
\begin{equation*}
	g_{\rm ess} = g \mathds{1}_{(\vre, \vte) \in K},\ 
	g_{\rm res} = g - g_{\rm ess} = g \mathds{1}_{(\vre, \vte) \in \R^2 \setminus K}
\end{equation*}
for any measurable $g = g(t,x)$. This decomposition obviously depends on $\ep$. The characteristic function $\mathds{1}_{(\vre, \vte) \in K}$ can be replaced by its smooth regularization by a suitable convolution kernel.

Here, we consider
\[
K = \Ov{\mathcal{U}(\Ov{\vr}, \Ov{\vt})} \subset (0, \infty)^2,\ { \mathcal{U} (\Ov{\vr}, \Ov{\vt} )} \ \mbox{- an open neighborhood of}\ (\Ov{\vr}, \Ov{\vt}).
\]
As shown in \cite[Chapter 5, Lemma~5.1]{FeNo6A}, there is a positive constant $C$ such that
\begin{equation} \label{u5}
E_{\ep} \left( \vr, \vt, \vu \Big| \tvr, \tvt, \tvu \right) \geq 
C \left( \frac{ |\vr - \tvr|^2 }{\ep^2} + \frac{ |\vt - \tvt|^2 }{\ep^2} + |\vu - \tvu |^2 \right)
\end{equation}
if $(\vr, \vt) \in K = \Ov{\mathcal{U}(\Ov{\vr}, \Ov{\vt})}$, $(\tvr, \tvt) \in \mathcal{U}(\Ov{\vr}, \Ov{\vt})$, and
\begin{equation} \label{u6}
E_{\ep} \left( \vr, \vt, \vu \Big| \tvr, \tvt, \tvu \right) \geq 
C \left( \frac{1}{\ep^2} + \frac{1}{\ep^2} \vr e(\vr, \vt) + \frac{1}{\ep^2} \vr |s(\vr, \vt)| + \vr |\vu|^2 \right)
\end{equation}
whenever $(\vr, \vt) \in R^2 \setminus \Ov{\mathcal{U}(\Ov{\vr}, \Ov{\vt})}$, $(\tvr, \tvt) \in \mathcal{U}(\Ov{\vr}, \Ov{\vt})$. The constant $C$
depends on the compact set $K$ and the distance 
\[
\sup_{t,x} {\rm dist} \left[ (\tvr (t,x), \tvt (t,x) ) ; \partial K \right]. 	\]

\subsection{Energy estimates}

The necessary energy bounds are obtained by plugging 
\[
(\tvr, \tvt, \tvu) = (\Ov{\vr}, \Ov{\vt} + \ep \vtB, 0)
\]
as ``test'' functions in the relative energy inequality \eqref{u3}. Here $\vtB = \vtB(x)$ is the $C^2$ function generating the temperature 
boundary data, keep in mind \eqref{eq:t_boundary}. After a straightforward manipulation, we obtain  
 \begin{align}
	&\left[ \intO{ E_\ep \left(\vre, \vte, \vue \Big| \Ov{\vr}, \Ov{\vt} + \ep \vtB, 0 \right) } \right]_{t = 0}^{t = \tau} \br 
	&+ \int_0^\tau  \int_{{\Omega}} \frac{\Ov{\vt} + \ep \vtB }{\vte} \left(  \mathbb{S}(\vte, \Ds \vue) : \Ds \vue + \frac{1}{\ep^2} \frac{\kappa(\vte) |\Grad \vte|^2 }{\vte}\right)   \dx \dt \br 
	&\leq - \frac{1}{\ep} \int_0^\tau \intO{ \left( \vre (s(\vre,\vte) - s(\Ov{\vr}, \Ov{\vt} + \ep \vtB)) \vue \cdot \Grad \vtB \right) } \dt + \frac{1}{\ep} \int_0^\tau \intO{ 
		\frac{\kappa (\vte) \Grad \vte}{\vte}  \cdot \Grad \vtB } \dt \br 
	&+ \int_0^\tau \intO{ \vre \left[ \frac{1}{2\ep} \Grad G + \frac{1}{\ep} \Grad |\vc{x}_h |^2 
		- \frac{2}{\sqrt{\ep}} (\vc{e}_3 \times \vue)   \right] \cdot \vue  } \dt  \br 
	&- \frac{1}{\ep} \int_0^\tau \intO{ \left[  \left( \frac{\partial p(\Ov{\vr}, \Ov{\vt} + \ep \vtB) }{\partial \vt} - 
		\frac{\partial p(\Ov{\vr}, \Ov{\vt})}{\partial \vt} \right)\frac{\vre}{\tvr} \vue \cdot \Grad \vtB \right] } \dt \br 
	&- \frac{1}{\ep} \int_0^\tau \intO{ \frac{\partial p(\Ov{\vr}, \Ov{\vt})}{\partial \vt} \frac{\vre}{\tvr} \vue \cdot \Grad \vtB} \dt.	
	\label{ee1}
\end{align}

Thanks to hypothesis \eqref{LLw1}, 
\[
\intO{ E_\ep \left(\vr_{0,\ep}, \vt_{0,\ep}, \vu_{0,\ep} \Big| \Ov{\vr}, \Ov{\vt} + \ep \vtB, 0 \right) } \aleq 1.
\]
Moreover, since 
\[
\frac{1}{\sqrt{\ep}} (\vc{e}_3 \times \vue) \cdot \vue = 0,  
\]
inequality \eqref{ee1} coincides with its counterpart in \cite[Section 5, formula 5.4]{BelFeiOsch}. 

Thanks to our choice 
of the velocity boundary conditions \eqref{B1}, \eqref{B2}, and hypothesis \eqref{w16}, we have, similarly to \cite{BelFeiOsch},
the Korn-Poincar\' e inequality, 
\begin{equation} \label{KP}
\| \vu \|^2_{W^{1,2} (\Omega; R^3)} \aleq \intO{ \frac{1}{\vt} \mathbb{S}(\vt, \Ds \vu) : \Ds \vu }. 
\end{equation}	
Consequently, we may repeat step by step the arguments of \cite[Section 5.1]{BelFeiOsch} to obtain the following list 
of uniform bounds, cf. \cite[Section 5.1.2]{BelFeiOsch}: 
\begin{align}
	{\rm ess} \sup_{t \in (0,T)} \intO{ E_\ep \left( \vre, \vte, \vue \Big| \Ov{\vr}, \Ov{\vt} + \ep \vtB, 0 \right) } &\aleq 1, \label{UB1} \\
	\int_0^T \| \vue \|^2_{W^{1,2} (\Omega; \mathbb{R}^3) } \dt &\aleq 1, \label{UB2} \\
	\frac{1}{\ep^2} \int_0^T \left( \| \Grad \log (\vte) \|^2_{L^2(\Omega; \mathbb{R}^3)} + \| \Grad \vte^{\frac{\beta}{2}} \|^2_{L^2(\Omega; \mathbb{R}^3)} \right) &\aleq 1
	\label{UB3}
\end{align}
uniformly for $\ep \to 0$. Moreover, using the structural hypotheses imposed on the EOS we deduce,  
\begin{align} \
	\frac{1}{\ep^2} {\rm ess} \sup_{t \in (0,T)} \intO{ [1]_{\rm res} } &\aleq 1 \label{UB4} \\ 
	{\rm ess} \sup_{t \in (0,T)} \intO{ \vre |\vue|^2 } &\aleq 1, \label{UB5} \\ 
	{\rm ess} \sup_{t \in (0,T)} \left\| \left[ \frac{\vre - \Ov{\vr}}{\ep} \right]_{\rm ess} \right\|_{L^2(\Omega)} &\aleq 1, \label{UB6} 
	\\
	{\rm ess} \sup_{t \in (0,T)} \left\| \left[ \frac{\vte - \Ov{\vt}}{\ep} \right]_{\rm ess} \right\|_{L^2(\Omega)} &\aleq 1, \label{UB7} 
	\\
	\frac{1}{\ep^2} {\rm ess} \sup_{t \in (0,T)} \| [\vre]_{\rm res} \|^{\frac{5}{3}}_{L^{\frac{5}{3}}(\Omega)} + \frac{1}{\ep^2}  {\rm ess} \sup_{t \in (0,T)} \| [\vte]_{\rm res} \|^{4}_{L^{4}(\Omega)}	 &\aleq 1.
	\label{UB8}	
\end{align}
Finally, combining the above bounds, we have 
\begin{equation} \label{UB9}
	\int_0^T \left\| \frac{\log(\vte) - \log(\Ov{\vt})}{\ep} \right\|^2_{W^{1,2}(\Omega)} \dt + \int_0^T \left\| \frac{ \vte - \Ov{\vt} }{\ep} \right\|^2_{W^{1,2}(\Omega)} \dt \aleq 1,
\end{equation}
and 
\begin{equation} \label{UB10} 
	\int_0^T \left\| \left[ \frac{\kappa (\vte) }{\vte} \right]_{\rm res} \frac{\Grad \vte }{\ep} \right\|^q_{L^q(\Omega; \mathbb{R}^d)} \dt \aleq 1 \ \mbox{for some}\ q > 1, 
\end{equation} 
\begin{equation} \label{UB11}
	\| \vte^{\frac{\beta}{2}} \|_{L^r ((0,T) \times \Omega)} \aleq 1 \ \mbox{for some}\ r > 2
\end{equation}
uniformly for $\ep \to 0$, cf. \cite[Section 5.1.2]{BelFeiOsch}.

\section{Asymptotic limit}
\label{A}

With the uniform bounds established in the preceding section, we are able to perform the limit 
$\ep \to 0$ and to identify the target system.

In the first section, we will identify (up to a suitable extraction) limit points (with respect to suitable weak topologies)
of the sequence $\big(\vre, \vte, \vue\big)_\ep$ of weak solutions of the NSF system. More precisely,
we will prove weak convergence properties for the sequence of the velocity fields $\vue$, and of the quantities $(\vre - \Ov{\vr})/\ep$ and $(\vte - \Ov{\vt})/\ep$.
The limit points of those families will be denoted, respectively, by $\vu$, $\mathfrak{R}$ and $\mathfrak{T}$. 
Recall that the solution of the target problem \eqref{TS1}--\eqref{TS5} is instead denoted by $\big(\mathcal{R}, \mathcal{T}, \vU\big)$,
with $\vU = (\vuh, 0)$. See also \eqref{A12} and \eqref{A13} below.

Notice that, at this point, there is no reason why the triplets $\big(\mathfrak{R}, \mathfrak{T}, \vu\big)$ and $\big(\mathcal{R}, \mathcal{T}, \vU\big)$
should coincide. This will be a consequence of the computations of Section \ref{s:strong}, where, by means of a relative energy argument,
we will show strong convergence of  $\big(\frac{\vre - \Ov{\vr}}{\ep}, \frac{\vte - \Ov{\vt}}{\ep} , \vue\big)_\ep$ to
$\big(\mathcal{R}, \mathcal{T}, \vU\big)$.
In particular, this will yield the proof of Theorem \ref{TM1}.

\subsection{Weak convergence}

First, it follows from \eqref{UB2}--\eqref{UB8} that 
\begin{align} 
	\vre &\to \Ov{\vr} \ \mbox{in}\ L^{\frac{5}{3}}(\Omega) \ \mbox{uniformly for}\ t \in (0,T),  \label{A1} \\ 
	\vte &\to \Ov{\vt} \ \mbox{in}\ L^2(0,T; W^{1,2}(\Omega)), \label {A2} \\
	\vue &\to \vu \ \mbox{weakly in}\ L^2(0,T; W^{1,2}(\Omega; \mathbb{R}^3)), \label{A3}
\end{align}
where $\vu$ satisfies the boundary conditions 
\begin{equation} \label{A4} 
\vu|_{\partial B(r) \times (0,1)} = 0,\ 
\vu \cdot \vc{n}|_{x_3 =1,2} = 0.	
\end{equation}	
Strictly speaking, the limit \eqref{A3} holds modulo a suitable subsequence.
Owing to Proposition \ref{TP1},
however, the limit $\vu = \vU$ is uniquely determined by the target system and the convergence is therefore unconditional.

Next, by the same token, we may extract a suitable subsequence so that 
\begin{align}
	\frac{ \vre - \Ov{\vr} }{\ep} &= \left[ \frac{ \vre - \Ov{\vr} }{\ep} \right]_{\rm ess} + 
	\left[ \frac{ \vre - \Ov{\vr} }{\ep} \right]_{\rm res}, \br 
	\left[ \frac{ \vre - \Ov{\vr} }{\ep} \right]_{\rm ess} &\to \mathfrak{R} \ \mbox{weakly-(*) in}\ 
	L^\infty(0,T; L^2(\Omega)),\ \intO{ \mathfrak{R}(t, \cdot) } = 0, \br
	\left[ \frac{ \vre - \Ov{\vr} }{\ep} \right]_{\rm res} &\to 0 \ \mbox{in}\ L^\infty(0,T; L^{\frac{5}{3}}(\Omega)),
	\label{A5}	
\end{align}
and 
\begin{equation}
	\frac{ \vte - \Ov{\vt} }{\ep} \to  \mathfrak{T} \ \mbox{weakly in}\  L^2(0,T; W^{1,2}(\Omega)) 
	\ \mbox{and weakly-(*) in}\ L^\infty(0,T; L^2(\Omega)),
	\label{A6}	
\end{equation}
where
\begin{equation} \label{A7}
	\mathfrak{T}|_{\partial \Omega} = \vtB.
\end{equation}

Finally, we perform the limit in the equation of continuity \eqref{Lw4} obtaining 
\begin{equation} \label{A8}
\Div \vu = 0, 
\end{equation}
while the limit in the momentum equation yields 
\begin{equation} \label{A9}
	\frac{\partial p(\Ov{\vr}, \Ov{\vt})}{\partial \vr} \Grad \mathfrak{R} + \frac{\partial p(\Ov{\vr}, \Ov{\vt})}{\partial \vt} \Grad \mathfrak{T} =
	\Ov{\vr} \Big( \Grad G + \frac{1}{2} \Grad |\mbf{x}_h|^2 \Big).
\end{equation}	
In particular, we deduce that
\begin{equation} \label{A10}
	\mathfrak{R} \in L^2(0,T; W^{1,2}(\Omega)).
\end{equation}

\subsubsection{Target velocity profile}

Now, we test momentum equation \eqref{Lw5} on $\sqrt{\ep} \Curl \bfphi$, $\bfphi \in C^1_c((0,T) \times \Omega;R^3)$. Seeing that 
\[
\frac{1}{\sqrt{\ep}} \intO{ \vre \Big( \Grad G + \frac{1}{2}\Grad |\vc x_h|^2 \Big) \Curl \bfphi } = 
\sqrt{\ep} \intO{ \frac{\vre - \Ov{\vr} }{\ep} \Big( \Grad G + \frac{1}{2} \Grad |\vc x_h|^2 \Big) \Curl \bfphi } \to 0
\]
as $\ep \to 0$, we conclude
\[
\int_0^T \intO{ \Ov{\vr} (\vc e_3 \times \vu) \cdot \Curl \bfphi } \dt = 0.  
\]	
Thus 
\[
\Curl (\vc e_3 \times \vu) = 0 \ \Rightarrow \ [- u^2, u^1, 0] = \Grad \Psi 
\]
on any simply connected subset of $\Omega$. This yields $\Psi$ is independent of $x_3$; whence 
\[
[u^1, u^2] = [u^1, u^2] (t, x_h),\ \Divh (u^1, u^2) = 0 \ \Rightarrow \ u^3 \ \mbox{independent of}\ x_3,
\]
where the last implication follows from \eqref{A8}.
However, in accordance with \eqref{A4}, $u^3$ vanishes for $x_3 = 0,1$; whence $u^3 = 0$. We therefore conclude 
that the limit velocity profile satisfies 
\begin{equation} \label{A11}
\vu = (u^1, u^2, 0),\ (u^1, u^2) = 	(u^1, u^2) (t,\vc x_h),\ \Divh (u^1, u^2) = 0.
\end{equation}	

Thus we have observed that fast rotation forces the limit velocity to be two-dimensional and purely horizontal, in the sense that it depends
only on the horizontal variable $\vc{x}_h$. This can be 
seen as a mathematical formulation of the celebrated \emph{Taylor-Proudman theorem}
in geophysics \cite{GILL, PEDL}.

\subsection{Strong convergence to the target system} \label{s:strong}

To complete the proof of Theorem \ref{TM1}, it remains to establish the strong convergence to the target problem. This will be achieved 
by considering the limit as the test function in the relative energy inequality \eqref{u3}. 

First, it is convenient to rewrite the target system in terms of the variables $(\vuh, \mathcal{R}, \mathcal{T})$,  where 
\begin{equation} \label{A12}
\mathcal{T} = \Theta + \frac{\lambda(\Ov{\vr}, \Ov{\vt})} {1 - \lambda(\Ov{\vr}, \Ov{\vt})} \avintO{ \Theta } 
\end{equation}
Accordingly, the target problem reads
\begin{align} 
\Divh \vuh &= 0, \br	
\Ov{\vr} \Big[ \partial_t \vuh + \Divh (\vuh \otimes \vuh) \Big] + \Gradh \Pi &= 
\mu(\Ov{\vt}) \Delta_h \vuh + \left<\mathcal{R} \right> \Gradh \Big( G + \frac{1}{2}|\vc x_h|^2 \Big), \br 
\Ov{\vr} c_p(\Ov{\vr}, \Ov{\vt} ) \Big[ \partial_t \mathcal{T} + \vuh \cdot \Gradh \mathcal{T} \Big] &- 
\Ov{\vr} \Ov{\vt} \alpha (\Ov{\vr}, \Ov{\vt} ) \vuh \cdot \Gradh \Big( G + \frac{1}{2}|\vc x_h|^2 \Big) \br &= 
\kappa (\Ov{\vt}) \Del \mathcal{T} + \Ov{\vt}  \alpha (\Ov{\vr}, \Ov{\vt} )  \frac{\partial p  (\Ov{\vr}, \Ov{\vt} ) } 
{\partial \vt} \partial_t \avintO{ \mathcal{T} },\ \br 
	\frac{\partial p(\Ov{\vr}, \Ov{\vt})}{\partial \vr} \Grad \mathcal{R} + \frac{\partial p(\Ov{\vr}, \Ov{\vt})}{\partial \vt} \Grad \mathcal{T} &=
\Ov{\vr} \Big( \Grad G + \frac{1}{2} \Grad |\vc{x}_h|^2 \Big),\ \intO{ \mathcal{R} } = 0,
\label{A13}	
\end{align}	
with the boundary conditions
\begin{equation} \label{A14}
\vuh|_{\partial B(r)} = 0,\ 
\mathcal{T}|_{\partial \Omega} = \vtB, 
\end{equation}

\subsubsection{Estimates based on the relative energy inequality}

To finish the proof of Theorem \ref{TM1}, we consider the trio 
\[
\tvr = \Ov{\vr} + \ep \mathcal{R},\ \tvt = \Ov{\vt} + \ep \mathcal{T}, \ \tvu = \vU = (\vuh, 0) 
\]
as a ``test function'' in the relative energy inequality \eqref{u3}. As the initial data are well--prepared, we get
\begin{equation} \label{A15}
	\intO{ E_\ep \left(\vre, \vte, \vue \ \Big| \Ov{\vr} + \ep \mathcal{R}, \Ov{\vt} + \ep \mathcal{T}, \vU \right) (0, \cdot) 
	} \to 0 \ \mbox{as}\ \ep \to 0.
\end{equation}

The next observation is that the integral corresponding to the Coriolis force vanishes. Indeed, owing to the fact that 
$\vU = (\vuh, 0)$, with $\vuh = \vuh(t,\vc x_h)$ and $\Divh \vuh=0$, there exists
a stream-function $\Phi = \Phi(t,\vc x_h)$ such that
\[
 \vuh = \nabla_h^\perp\Phi = (-\partial_2\Phi, \partial_1 \Phi, 0).
\]
However, noticing that $\vc e_3\times \vU = ( -u_h^2, u_h^1, 0 ) = - (\partial_1\Phi, \partial_2\Phi, 0)$ and that $\partial_3\Phi=0$, we can write
\begin{equation} \label{eq:potential}
 \vc e_3\times \vU = - \Grad\Phi, \qquad \mbox{ with } \ \Phi = \Phi(t,\vc x_h).
\end{equation}
Thus, we may
use the weak formulation of the equation of continuity to write
\begin{align}
\int_0^\tau &\intO{ \frac{\vre}{\sqrt{\ep}} (\vc e_3 \times \vue) \cdot (\vue - \vU) } \dt = 
\int_0^\tau \intO{ \frac{\vre}{\sqrt{\ep}} \vc e_3\times \vU \cdot \vue } \br & = 
-\frac{1}{\sqrt{\ep}} \int_0^\tau \intO{ \vre \vue \cdot \Grad \Phi } \dt \br 
& = - \sqrt{\ep} \int_0^\tau \intO{ \frac{\vre - \Ov{\vr}}{\ep}  \partial_t \Phi } \dt + \sqrt{\ep} \left[ \intO{ \frac{\vre - \Ov{\vr} }{\ep} \Phi } \right]_{t = 0}^{t = \tau}
\label{A16}
\end{align}
where the right--hand side vanishes for $\ep \to 0$ as a consequence of the uniform bounds \eqref{UB4}, \eqref{UB6}, \eqref{UB8}. 
%Here, the symbol $\Phi$ denotes the potential associated to the irrotational field $(- u^2_h, u^1_h, 0)$.

\medskip

\noindent
$\bullet$ In view of the above observations, the relative energy inequality \eqref{u3} takes the form 
 \begin{align}
	&\intO{ E_\ep \left(\vre, \vte, \vue \Big| \Ov{\vr} + \ep \mathcal{R} , \Ov{\vt} + \ep \mathcal{T}, \vU \right) (\tau, \cdot) }  \br 
	&+ \int_0^\tau  \int_{{\Omega}} \frac{\Ov{\vt} + \ep \mathcal{T} }{\vte} \left(  \mathbb{S}(\vte, \Ds \vue) : \Ds \vue + \frac{1}{\ep^2} \frac{\kappa (\vte) \Grad \vte \cdot \Grad \vte}{\vte}\right)   \dx \dt \br 
	&\leq - \frac{1}{\ep} \int_0^\tau \intO{ \vre \Big( s(\vre, \vte) - s(\Ov{\vr} + \ep \mathcal{R}, \Ov{\vt} + \ep \mathcal{T}) \Big) \partial_t \mathcal{T} } \dt \br 
	&\quad - \frac{1}{\ep} \int_0^\tau \intO{ \vre \Big( s(\vre, \vte) - s(\Ov{\vr} + \ep \mathcal{R}, \Ov{\vt} + \ep \mathcal{T}) \Big) \vue \cdot \Grad \mathcal{T}  } \dt
	 \br 
	&\quad + \frac{1}{\ep} \int_0^\tau \intO{ 
		\frac{\kappa (\vte) \Grad \vte}{\vte}  \cdot \Grad \mathcal{T} } \dt \br 
	&\quad - \int_0^\tau \intO{ \Big[ \vre (\vue - \vU) \otimes (\vue - \vU)  - \mathbb{S}(\vte, \Ds \vue) \Big] : \Ds \vU } \dt   \br
	&\quad + \int_0^\tau \intO{ \vre \left[ \frac{1}{\ep} \Grad G + \frac{1}{2\ep} \Grad |\vc{x}_h |^2 
	  - \partial_t \vU - (\vU \cdot \Grad) \vU  \right] \cdot (\vue - \vU) } \dt  \br 
	&+ \frac{1}{\ep^2} \int_0^\tau \intO{ \left[ \left( 1 - \frac{\vre}{\Ov{\vr} + \ep \mathcal{R}} \right) 
	\partial_t p(\Ov{\vr} + \ep \mathcal{R}, \Ov{\vt} + \ep \mathcal{T}) - \frac{\vre}{\Ov{\vr} + \ep \mathcal{R}} \vue \cdot \Grad p(\Ov{\vr} + \ep \mathcal{R}, \Ov{\vt} + \ep \mathcal{T}) \right] } \dt + \omega_\ep(\tau) ,
	\label{A17}
\end{align}
where the time--dependent function $\o_\ep(\tau)$ satisfies
\begin{equation} \label{def:omega_ep}
\sup_{\tau\in[0,T]}\omega_\ep(\tau) \to 0 \ \mbox{as}\ \ep \to 0, \quad \mbox{ for any given } T>0.
\end{equation}

\medskip 

\noindent
$\bullet$ Next, recalling that $\vU = (\vuh,0)$, we can use the limit momentum equation to rewrite 
\begin{align}
\int_0^\tau &\intO{ \vre \left[ \frac{1}{\ep} \Grad G + \frac{1}{2\ep} \Grad |\vc{x}_h |^2 
	- \partial_t \vU - (\vU \cdot \Grad) \vU  \right] \cdot (\vue - \vU) } \dt \br
&= \int_0^\tau \intO{ \frac{\vre}{\Ov{\vr}} \left[ \frac{1}{\ep} \Ov{\vr} \Big( \Grad G +  \frac{1}{2} \Grad |\vc{x}_h |^2 \Big) 
	- \Ov{\vr} \Big( \partial_t \vU - (\vU \cdot \Grad) \vU \Big)  \right] \cdot (\vue - \vU) } \dt	\br 
&=	\int_0^\tau \intO{ \frac{\vre}{\Ov{\vr}} \frac{1}{\ep} \Ov{\vr} \Big( \Grad G +  \frac{1}{2}\Grad |\vc{x}_h |^2 \Big) \cdot (\vue - \vU) } \dt \br
&+ \int_0^\tau \intO{ \frac{\vre}{\Ov{\vr}} \Big( \Gradh \Pi - \mu(\Ov{\vt}) \Delta_h \vuh - 
\left< \mathcal{R} \right> \Big( \nabla_h G +  \frac{1}{2} \nabla_h |\vc{x}_h |^2 \Big) \Big)   \cdot (\vu_{\ep,h} - \vuh) } \dt.
\nonumber
\end{align}
In view of the available uniform bounds, we may pass to the limit in the second integral, obtaining 
\begin{align} 
\int_0^\tau &\intO{ \frac{\vre}{\Ov{\vr}} \Big( \Gradh \Pi - \mu(\Ov{\vt}) \Delta_h \vuh - 
\left< \mathcal{R} \right> \Big( \nabla_h G +  \frac{1}{2} \nabla_h |\vc{x}_h |^2 \Big) \Big)   \cdot (\vu_{\ep,h} - \vuh) } \dt
%\int_0^\tau &\intO{ \frac{\vre}{\Ov{\vr}} \Big( \Gradh \Pi - \mu(\Ov{\vt}) \Delta_h \vuh - \left< r \right> \Big( \Grad G +  \frac {1}{2} \Grad |\vc{x}_h |^2 \Big) \Big)   \cdot (\vue - \vuh) } \dt 
\br
&=  \int_0^\tau \int_{B(r)}{ \Big( \Gradh \Pi - \mu(\Ov{\vt}) \Delta_h \vuh - 
\left< \mathcal{R} \right> \Big( \nabla_h G +  \frac{1}{2} \nabla_h |\vc{x}_h |^2 \Big) \Big)   \cdot \lan \vu_{\ep,h} - \vuh \ran } {\rm d}\vc x_h \dt
%\int_0^\tau \intO{  \Big( \Gradh \Pi - \mu(\Ov{\vt}) \Delta_h \vuh - \left< r \right> \Big( \Grad G +  \frac{1}{2}\Grad |\vc{x}_h |^2 \Big) \Big)   \cdot (\vu - \vuh) } \dt
+ \omega_\ep(\tau) \br 
&= %\int_0^\tau \intO{  \Big( - \mu(\Ov{\vt}) \Delta_h \vuh -  r  \Big( \Grad G +  \frac{1}{2} \Grad |\vc{x}_h |^2 \Big) \Big)   \cdot (\vu - \vuh) } \dt + \omega(\ep),
 \int_0^\tau \int_{B(r)}{ \Big( - \mu(\Ov{\vt}) \Delta_h \vuh -
\left< \mathcal{R} \right> \Big( \nabla_h G +  \frac{1}{2} \nabla_h |\vc{x}_h |^2 \Big) \Big)   \cdot \lan \vu_{\ep,h} - \vuh \ran } {\rm d}\vc x_h \dt 
+ \omega_\ep(\tau), 
\nonumber	
\end{align}	
where, similarly to \cite[Section 2, Proposition 2.1]{Ch-F-Gall}, we have used that
\[
\int_0^\tau \int_{B(r)}{  \Gradh \Pi   \cdot \lan \vu_{\ep,h} - \vuh \ran } {\rm d}\vc x_h \dt  = \omega_\ep(\tau),
\]
in the sense of \eqref{def:omega_ep}. Next, we observe that, since both $\Grad G$ and $\vU = ( \vuh , 0 )$ are independent of the vertical variable $x_3$,
the following equality holds:
\begin{align*}
&\int_0^\tau \int_{B(r)}{ \Big( - \mu(\Ov{\vt}) \Delta_h \vuh -
\left< \mathcal{R} \right> \Big( \nabla_h G +  \frac{1}{2} \nabla_h |\vc{x}_h |^2 \Big) \Big)   \cdot \lan \vu_{\ep,h} - \vuh \ran } {\rm d}\vc x_h \dt \br
&\qquad\qquad
= \int_0^\tau \intO{ \Big( - \mu(\Ov{\vt}) \Delta \vU -
\left< \mathcal{R} \right> \Big( \Grad G +  \frac{1}{2} \Grad |\vc{x}_h |^2 \Big) \Big)   \cdot ( \vu_{\ep} - \vU ) } \dt
\end{align*}

Consequently, inequality \eqref{A17} reduces to 
 \begin{align}
	&\intO{ E_\ep \left(\vre, \vte, \vue \Big| \Ov{\vr} + \ep \mathcal{R} , \Ov{\vt} + \ep \mathcal{T}, \vU \right) (\tau, \cdot) }  \br
	&+ \int_0^\tau \intO{ \Big( \mathbb{S}(\Ov{\vt}, \Ds \vue) - \mathbb{S}(\Ov{\vt}, 
		\Ds \vU) \Big): \Big( \Ds \vue - \Ds \vU \Big)    } \dt \br
	&+ \int_0^\tau  \int_{{\Omega}} \left( \frac{\Ov{\vt} + \ep \mathcal{T} }{\vte^2} \right)  \frac{\kappa (\vte) \Grad \vte \cdot \Grad \vte}{\ep^2 }   \dx \dt - \frac{1}{\ep} \int_0^\tau \intO{ 
		\frac{\kappa (\vte) \Grad \vte}{\vte}  \cdot \Grad \mathcal{T} } \dt\br 
	&\leq - \frac{1}{\ep} \int_0^\tau \intO{ \vre \Big( s(\vre, \vte) - s(\Ov{\vr} + \ep \mathcal{R}, \Ov{\vt} + \ep \mathcal{T}) \Big) \partial_t \mathcal{T} } \dt \br 
	&\quad - \frac{1}{\ep} \int_0^\tau \intO{ \vre \Big( s(\vre, \vte) - s(\Ov{\vr} + \ep \mathcal{R}, \Ov{\vt} + \ep \mathcal{T}) \Big) \vue \cdot \Grad \mathcal{T}  } \dt
	\br 
	&+ \int_0^\tau \intO{ \frac{\vre}{\Ov{\vr}} \frac{1}{\ep} \Ov{\vr} \Big( \Grad G + \frac{1}{2}\Grad |\vc{x}_h |^2 \Big) \Big( \vue - \vU \Big) } \dt \br
    &- \int_0^\tau \intO{ \lan\mathcal{R}\ran \Big( \Grad G + \frac{1}{2} \Grad |\vc{x}_h |^2 \Big) \cdot \Big( \vue - \vU \Big) } \dt  \br 
	&+ \frac{1}{\ep^2} \int_0^\tau \intO{ \left[ \left( 1 - \frac{\vre}{\Ov{\vr} + \ep \mathcal{R}} \right) \partial_t p(\Ov{\vr} + \ep \mathcal{R}, \Ov{\vt} + \ep \mathcal{T}) - \frac{\vre}{\Ov{\vr} + \ep \mathcal{R}} \vue \cdot \Grad p(\Ov{\vr} + \ep \mathcal{R}, \Ov{\vt} + \ep \mathcal{T}) \right] } \dt \br &+ C \int_0^\tau \intO{ E_\ep \left(\vre, \vte, \vue \Big| \Ov{\vr} + \ep \mathcal{R} , \Ov{\vt} + \ep \mathcal{T}, \vU \right) (\tau, \cdot) } \dt + 
	\omega_\ep(\tau).
	\label{A18}
\end{align}

\medskip 

\noindent
$\bullet$ The next step is integrating the Boussinesq relation \eqref{TS3}, obtaining
\[
\frac{\partial p(\Ov{\vr}, \Ov{\vt} ) }{\partial \vr}  \mathcal{R} + 
\frac{\partial p(\Ov{\vr}, \Ov{\vt} ) }{\partial \vt}  \mathcal{T} = \Ov{\vr} \Big( G + \frac{1}{2}|\vc x_h|^2 \Big) + \chi(t), 
\]
where we can replace $G$ by $G + {\rm const}$ and suppose that
\[
\intO{ \Big( G + \frac{1}{2}|\vc x_h|^2 \Big) } = 0.
\]
As $\mathcal{R}$ has zero average,
\[
\intO{ \mathcal{R} } = 0, 
\]
we get that
\[
\chi(t) = \frac{\partial p(\Ov{\vr}, \Ov{\vt})}{\partial \vt} \avintO{ \mathcal{T} (t, \cdot) }.
\]
Consequently, we may compute 
\begin{align}
	\frac{1}{\ep^2} & \intO{ \left( 1 - \frac{\vre}{\Ov{\vr} + \ep \mathcal{R}} \right) \partial_t p (\Ov{\vr} + \ep \mathcal{R}, \Ov{\vt} + \ep \MTC) } \br &= 
	\frac{1}{\ep} \intO{ \left( 1 - \frac{\vre}{\Ov{\vr} + \ep \mathcal{R}} \right) \left( \frac{\partial p (\Ov{\vr} + \ep \mathcal{R}, \Ov{\vt} + \ep \MTC) }{\partial \vr} \partial_t \mathcal{R} +  \frac{\partial p (\Ov{\vr} + \ep \mathcal{R}, \Ov{\vt} + \ep \MTC) }{\partial \vt} \partial_t \MTC  \right)	} \br 
	&= \intO{ \frac{1}{\ep} \left( 1 - \frac{\vre}{\Ov{\vr} + \ep \mathcal{R}} \right)  \left( \frac{\partial p(\Ov{\vr} + \ep \mathcal{R}, \Ov{\vt} + \ep \MTC)}{\partial \vr} -\frac{\partial p (\Ov{\vr} , \Ov{\vt} ) }{\partial \vr} \right)  \partial_t \mathcal{R} 
	} \br 
	&\quad + \intO{ \frac{1}{\ep}  \left( 1 - \frac{\vre}{\Ov{\vr} + \ep \mathcal{R}} \right) 
		\left( \frac{\partial p (\Ov{\vr} + \ep \mathcal{R}, \Ov{\vt} + \ep \MTC) }{\partial \vt} -\frac{\partial p (\Ov{\vr} , \Ov{\vt} ) }{\partial \vt} \right) \partial_t \MTC  	} \br
	&\quad + { \frac{1}{\ep} \intO{ \left( 1- \frac{\vre}{\Ov{\vr} + \ep \mathcal{R}} \right) \partial_t \chi }}.
	\label{A20}
\end{align}
Thanks to Proposition \ref{TP1}, straightforward computations show that
\begin{align}
	&\int^\tau_0 \intO{ \frac{1}{\ep} \left( 1 - \frac{\vre}{\Ov{\vr} + \ep \mathcal{R}} \right)  \left( \frac{\partial p(\Ov{\vr} + \ep \mathcal{R}, \Ov{\vt} + \ep \MTC)}{\partial \vr} -\frac{\partial p (\Ov{\vr} , \Ov{\vt} ) }{\partial \vr} \right)  \partial_t \mathcal{R} 
	} \dt \br 
	&\quad + \int_0^\tau \intO{ \frac{1}{\ep}  \left( 1 - \frac{\vre}{\Ov{\vr} + \ep \mathcal{R}} \right) 
		\left( \frac{\partial p (\Ov{\vr} + \ep \mathcal{R}, \Ov{\vt} + \ep \MTC) }{\partial \vt} -\frac{\partial p (\Ov{\vr} , \Ov{\vt} ) }{\partial \vt} \right) \partial_t \MTC  	} \dt
	= \omega_\ep(\tau),
	\nonumber
\end{align}
in the sense of relation \eqref{def:omega_ep}.
In addition, observing that
\begin{align*}
 \frac{1}{\ep} \left( 1 - \frac{\vre}{\Ov{\vr} + \ep \mathcal{R}} \right) &=
\frac{1}{\ep} \frac{\Ov{\vr} + \ep \mathcal{R} - \vre}{\Ov{\vr} + \ep \mathcal{R}} \br
&= -  \frac{\vre - \Ov{\vr}}{\ep( \Ov{\vr} + \ep \mathcal{R} )} + 
	\frac{\mathcal{R}}{\Ov{\vr} + \ep \mathcal{R}} \  \to \  \frac{1}{\Ov{\vr}} (\mathcal{R} - \mathfrak{R} ) \ \mbox{ as }\ \ep \to 0,
\end{align*}
where the convergence is in the weak-$*$ topology of $L^\infty(0,T;L^\frac53(\Omega))$,
and that 
\[
\intO{ \mathcal{R} } = \intO{ \mathfrak{R} } = 0,
\]
arguing again as in \cite[Proposition 2.1]{Ch-F-Gall} we may infer that
\[
\frac{1}{\ep} \int_0^\tau \intO{ \left( 1- \frac{\vre}{\Ov{\vr} + \ep \mathcal{R}} \right) \partial_t \chi } \dt = \omega_\ep(\tau).
\]

In a similar way, we may compute 
\begin{align}
	&- \frac{1}{\ep^2}  \int_0^\tau \intO{ \frac{\vre}{\Ov{\vr} + \ep \mathcal{R}} \vue \cdot \Grad p(\Ov{\vr} + \ep \mathcal{R}, \Ov{\vt} + \ep \MTC) } \dt \br &= 
	- \frac{1}{\ep} \int_0^\tau \intO{ \frac{\vre}{\Ov{\vr} + \ep \mathcal{R}} \vue \cdot \left( \frac{\partial p (\Ov{\vr} + \ep \mathcal{R}, \Ov{\vt} + \ep \MTC) }{\partial \vr} \Grad \mathcal{R} +  \frac{\partial p (\Ov{\vr} + \ep \mathcal{R}, \Ov{\vt} + \ep \MTC) }{\partial \vt} \Grad \MTC  \right)	} \dt \br 
	&=- \int_0^\tau \intO{ \frac{1}{\ep} \frac{\vre}{\Ov{\vr} + \ep \mathcal{R}} \vue \cdot \Grad \mathcal{R} \left( \frac{p \partial (\Ov{\vr} + \ep \mathcal{R}, \Ov{\vt} + \ep \MTC) }{\partial \vr} -\frac{\partial p (\Ov{\vr} , \Ov{\vt} ) }{\partial \vr} \right) } \dt \br
	&\quad - \int_0^\tau \intO{ \frac{1}{\ep} \frac{\vre}{\Ov{\vr} + \ep \mathcal{R}} \vue \cdot \Grad \MTC \left( \frac{\partial p (\Ov{\vr} + \ep \mathcal{R}, \Ov{\vt} + \ep \MTC) }{\partial \vt} -\frac{\partial p (\Ov{\vr} , \Ov{\vt} ) }{\partial \vt} \right) } \dt \br
	&- \frac{1}{\ep} \int_0^\tau \intO{\frac{\vre}{\Ov{\vr} + \ep \mathcal{R}} \Ov{\vr} \vue \cdot \Grad \Big( G + \frac{1}{2}|\vc x_h|^2 \Big)  } \dt .
	\label{A21}
\end{align}
In view of the convergence stated in \eqref{A1}, \eqref{A3}, and of the properties \eqref{A11} on the target velocity field $\vu$,
we get
\begin{align} 
	\int_0^\tau &\intO{ \frac{1}{\ep} \frac{\vre}{\Ov{\vr} + \ep \mathcal{R}} \vue \cdot \Grad \mathcal{R} \left( \frac{\partial p(\Ov{\vr} + \ep \mathcal{R}, \Ov{\vt} + \ep \MTC) }{\partial \vr} -\frac{\partial p (\Ov{\vr} , \Ov{\vt} ) }{\partial \vr} \right) } \dt \br
	& + \int_0^\tau \intO{ \frac{1}{\ep} \frac{\vre}{\Ov{\vr} + \ep \mathcal{R}} \vue \cdot \Grad \MTC \left( \frac{\partial p (\Ov{\vr} + \ep \mathcal{R}, \Ov{\vt} + \ep \MTC) }{\partial \vt} -\frac{\partial p (\Ov{\vr} , \Ov{\vt} ) }{\partial \vt} \right) } \dt \br
	&= \int_0^\tau \intO{ \vue \cdot \left(  \frac{\partial^2 p (\Ov{\vr} , \Ov{\vt}) }{\partial^2 \vr} \mathcal{R} \Grad \mathcal{R} + \frac{\partial^2 p (\Ov{\vr} , \Ov{\vt}) }{\partial \vr \partial \vt }  \Grad (\mathcal{R} \MTC) + 
		\frac{\partial^2 p (\Ov{\vr} , \Ov{\vt}) }{\partial^2 \vt} { \MTC} \Grad \MTC		\right) } + \omega_\ep(\tau) \br 
	&= \omega_\ep(\tau),
	\nonumber
\end{align}
where again we have argued as in \cite[Proposition 2.1]{Ch-F-Gall}.

Summarizing, we may rewrite inequality \eqref{A18} in the form
 \begin{align}
	&\intO{ E_\ep \left(\vre, \vte, \vue \Big| \Ov{\vr} + \ep \mathcal{R}, \Ov{\vt} + \ep \MTC, \vU \right) (\tau, \cdot) }  \br 
	&+ \int_0^\tau \intO{ \Big( \mathbb{S} (\Ov{\vt}, \Ds \vue) - \mathbb{S} (\Ov{\vt}, \Ds \vU) \Big) : \Big( \Ds \vue - \Ds \vU  \Big) } \dt \br &
	+\int_0^\tau \intO{  \left(  \frac{\Ov{\vt} + \ep \MTC}{\vte^2} \right) \frac{\kappa (\vte) \Grad \vte \cdot \Grad \vte }{\ep^2}  } \dt
	- \int_0^\tau \intO{ \frac{\kappa (\vte) }{\vte}  \frac{\Grad \vte }{\ep} \cdot \Grad \MTC      } \dt \br	
	&\leq - \frac{1}{\ep} \int_0^\tau \intO{ \vre \Big[ s(\vre, \vte)
		 - s(\Ov{\vr} + \ep \mathcal{R}, \Ov{\vt} + \ep \MTC) \Big] \partial_t \MTC } \dt \br  &\quad - \frac{1}{\ep} \int_0^\tau \intO{ \vre \Big[ s(\vre, \vte) - s(\Ov{\vr} + \ep \mathcal{R}, \Ov{\vt} + \ep \MTC) \Big] \vue \cdot \Grad \MTC   } \dt \br  
    &\quad + \int_0^\tau \intO{ \frac{\vre}{\Ov{\vr}} \frac{1}{\ep} \Ov{\vr} \Grad\Big( G + \frac{1}{2} |\vc{x}_h |^2 \Big)\cdot \vue } \dt
            - \int_0^\tau \intO{ \lan \mathcal{R} \ran \Grad \Big( G + \frac{1}{2} |\vc{x}_h |^2 \Big) \cdot \vue } \dt \br
    &\quad - \frac{1}{\ep} \int_0^\tau \intO{\frac{\vre}{\Ov{\vr} + \ep \mathcal{R}} \Ov{\vr} \vue \cdot \Grad \Big( G + \frac{1}{2}|\vc x_h|^2 \Big)  } \dt \br
    %&\quad - \int_0^\tau \intO{ \lan \mathcal{R} \ran  \Big( \Grad G + \frac{1}{2}\Grad |\vc{x}_h |^2 \Big) \vue } \dt \br
	&\quad - \int_0^\tau \intO{ \frac{\vre}{\ep}  \Grad \Big( G + \frac{1}{2}|\vc x_h|^2 \Big) \cdot \vU } \dt    +  
	\int_0^\tau \intO{ \lan \mathcal{R}\ran  \Grad \Big( G + \frac{1}{2}|\vc x_h|^2 \Big)  \cdot \vU } \dt  \br 
	&\quad + C \int_0^\tau \intO{E_\ep \left(\vre, \vte, \vue \Big| \Ov{\vr} + \ep r, \Ov{\vt} + \ep \MTC, \vU \right) } \dt   + \omega_\ep(\tau).
	\label{A22}
\end{align}

\medskip

\noindent
$\bullet$ Let us focus on the terms depending on $ \Grad\Big( G + \frac{1}{2}|\vc x_h|^2 \Big)$ in \eqref{A22}. First of all,
seeing that $ \vU = (\vuh, 0)$ is solenoidal, that is $ \Div \vU = \Divh \vuh = 0$, we may write
\begin{align}
-\int_0^\tau &\intO{ \frac{\vre}{\ep}  \Grad \Big( G + \frac{1}{2}|\vc x_h|^2 \Big) \cdot \vU } \dt    +  \int_0^\tau \intO{
\lan \mathcal{R}\ran \Grad \Big( G + \frac{1}{2}|\vc x_h|^2 \Big)  \cdot \vU } \dt \br
&=\int_0^\tau \intO{ \frac{ \vre - \Ov{\vr} - \ep \lan\mathcal{R}\ran}{\ep} \Grad \Big( G + \frac{1}{2}|\vc x_h|^2 \Big) \cdot \vU } \dt.
\label{eq:F-U}
\end{align}

In addition, we can compute
\begin{align*}
&\int_0^\tau \intO{ \frac{\vre}{\Ov{\vr}} \frac{1}{\ep} \Ov{\vr} \Big( \Grad G + \frac{1}{2}\Grad |\vc{x}_h |^2 \Big) \vue } \dt - \frac{1}{\ep} \int_0^\tau \intO{\frac{\vre}{\Ov{\vr} + \ep \mathcal{R}} \Ov{\vr} \vue \cdot \Grad \Big( G + \frac{1}{2}|\vc x_h|^2 \Big)  } \dt \br
&\qquad = \ \int_0^\tau \intO{ \frac{\vre \ \Ov{\vr}}{\ep} \left(\frac{1}{\Ov{\vr}} - \frac{1}{\Ov{\vr} + \ep \mathcal{R}}\right)
\Grad\Big( G + \frac{1}{2} |\vc{x}_h |^2 \Big)\cdot \vue } \dt \br
&\qquad = \ \int_0^\tau \intO{ \frac{\vre \ \mathcal{R}}{\Ov{\vr} + \ep \mathcal{R}} \Grad\Big( G + \frac{1}{2} |\vc{x}_h |^2 \Big)\cdot \vue } \dt.
\end{align*}
Now, owing to the convergence properties \eqref{A1}, \eqref{A3}, \eqref{A5}, we infer the series of equalities
\begin{align*}
\int_0^\tau \intO{ \frac{\vre \ \mathcal{R}}{\Ov{\vr} + \ep \mathcal{R}} \Grad\Big( G + \frac{1}{2} |\vc{x}_h |^2 \Big)\cdot \vue } \dt 
&= \int_0^\tau \intO{  \mathcal{R} \Grad\Big( G + \frac{1}{2} |\vc{x}_h |^2 \Big)\cdot \vue } \dt + \omega_\ep(\tau) \br
&= \int_0^\tau \intO{  \mathcal{R} \Grad\Big( G + \frac{1}{2} |\vc{x}_h |^2 \Big)\cdot \vu } \dt + \omega_\ep(\tau) \br
&= \int_0^\tau \intO{  \lan \mathcal{R} \ran \Grad\Big( G + \frac{1}{2} |\vc{x}_h |^2 \Big)\cdot \vu } \dt + \omega_\ep(\tau),
\end{align*}
where the last equality follows from the fact that $\vu = \vu (t,\vc  x_h)$, recall \eqref{A11}. We have thus proven that
\begin{align*}
&\int_0^\tau \intO{ \frac{\vre}{\Ov{\vr}} \frac{1}{\ep} \Ov{\vr} \Grad\Big( G + \frac{1}{2} |\vc{x}_h |^2 \Big)\cdot \vue } \dt
            - \int_0^\tau \intO{ \lan \mathcal{R} \ran \Grad \Big( G + \frac{1}{2} |\vc{x}_h |^2 \Big) \cdot \vue } \dt \br
    &\qquad\qquad\qquad\qquad\qquad\qquad\qquad
    - \frac{1}{\ep} \int_0^\tau \intO{\frac{\vre}{\Ov{\vr} + \ep \mathcal{R}} \Ov{\vr} \vue \cdot \Grad \Big( G + \frac{1}{2}|\vc x_h|^2 \Big)  } \dt 
    = \omega_\ep(\tau),
\end{align*}
in the sense of relation \eqref{def:omega_ep}.

Thus, combining equality \eqref{eq:F-U} and the previous relation with the estimates established in the preceding section and the convergence
\eqref{A5}, \eqref{A6}, we may rewrite  \eqref{A22} in the form 
 \begin{align}
	&\intO{ E_\ep \left(\vre, \vte, \vue \Big| \Ov{\vr} + \ep \mathcal{R}, \Ov{\vt} + \ep \MTC, \vU \right) (\tau, \cdot)}  \br 
	&+ \int_0^\tau \intO{ \Big( \mathbb{S} (\Ov{\vt}, \Ds \vue) - \mathbb{S} (\Ov{\vt}, \Ds \vU) \Big) : \Big( \Ds \vue - \Ds \vU  \Big) } \dt \br &
	+\int_0^\tau \intO{  \left(  \frac{\Ov{\vt} + \ep \MTC}{\vte^2} \right) \frac{\kappa (\vte) \Grad \vte \cdot \Grad \vte }{\ep^2}  } \dt
	- \int_0^\tau \intO{ \frac{\kappa (\Ov{\vt}) }{\Ov{\vt}} \Grad \mathfrak{T} \cdot \Grad { \MTC} } \dt \br	
	&\leq - \frac{1}{\ep} \int_0^\tau \intO{ \vre \Big[ (s(\vre, \vte) - s( \Ov{\vr} + \ep \mathcal{R},\Ov{\vt} + \ep \MTC ) \Big] \partial_t \MTC } \dt \br 
	&\quad - \frac{1}{\ep} \int_0^\tau \intO{ \vre \Big[ s(\vre, \vte) - s( \Ov{\vr} + \ep \mathcal{R}, \Ov{\vt} + \ep \MTC) \Big] \vue \cdot \Grad \MTC  } \dt \br
	&\quad  + \frac{1}{\ep} \int_0^\tau \intO{  \vre \Big[ s(\vre, \vte) - s(\Ov{\vr} + \ep \mathcal{R}, \Ov{\vt} + \ep \MTC) \Big] ( \vU - \vue) \cdot \Grad \MTC } \dt \br
	&\quad    +  \int_0^\tau \intO{ \lan \mathcal{R} - \mathfrak{R}\ran \Grad \Big( G + \frac{1}{2}|\vc x_h|^2 \Big)  \cdot \vU } \dt  \br 
	&\quad+ C \int_0^\tau \intO{E_\ep \left(\vre, \vte, \vue \Big| \Ov{\vr} + \ep \mathcal{R}, \Ov{\vt} + \ep \MTC, \vU \right) } \dt   + \omega_\ep(\tau).
	\label{A23}
\end{align}

\medskip

\noindent 
$\bullet$ Since the temperature deviation $\mathcal{T}$ satisfies the third equation in \eqref{A13}, we have 
\begin{align} 
	\partial_t \MTC + \vuh \cdot \Gradh \MTC &= \frac{\Ov{\vt} \alpha (\Ov{\vr}, \Ov{\vt} ) }{ c_p (\Ov{\vr}, \Ov{\vt} )} \Gradh \Big(G 
	+ \frac{1}{2}|\vc x_h|^2 \Big) \cdot \vuh + 
	\frac{\kappa(\Ov{\vt})}{\Ov{\vr} c_p (\Ov{\vr}, \Ov{\vt} )}  \Del \MTC 
	+ \frac{1}{\Ov{\vr} c_p (\Ov{\vr}, \Ov{\vt} )} \xi(t), \br 
	\ \mbox{where we have denoted }\ \xi(t)  &= 
	\Ov{\vt} \alpha (\Ov{\vr}, \Ov{\vt}) \frac{\partial p (\Ov{\vr}, \Ov{\vt})}{\partial \vt} 
	\partial_t \avintO{ \MTC(t,\cdot) }.	\label{A23a}
\end{align}
Thus, using the estimate
\begin{align}
	\frac{1}{\ep} &\int_0^\tau \intO{  \vre \Big[ s(\vr_\ep,\vt_\ep) - s(\Ov{\vr} + \ep \mathcal{R}, \Ov{\vt} + \ep \MTC) \Big] ( \vU - \vue) \cdot \Grad {\MTC} } \dt \br
	&\aleq \int_0^\tau \intO{E_\ep \left(\vre, \vte, \vue \Big| \Ov{\vr} + \ep \mathcal{R}, \Ov{\vt} + \ep \MTC, \vU \right) } \dt,
	\nonumber
\end{align}
we may rewrite \eqref{A23} in the form
 \begin{align}
	& \intO{ E_\ep \left(\vre, \vte, \vue \Big| \Ov{\vr} + \ep \mathcal{R}, \Ov{\vt} + \ep \MTC, \vU \right) (\tau, \cdot) }  \br 
	&+ \int_0^\tau \intO{ \Big( \mathbb{S} (\Ov{\vt}, \Ds \vue) - \mathbb{S} (\Ov{\vt}, \Ds \vU) \Big) : \Big( \Ds \vue - \Ds \vU  \Big) } \dt \br &
	+\int_0^\tau \intO{  \left(  \frac{\Ov{\vt} + \ep \MTC}{\vte^2} \right) \frac{\kappa (\vte) \Grad \vte \cdot \Grad \vte }{\ep^2}  } \dt
	- \int_0^\tau \intO{ \frac{\kappa (\Ov{\vt}) }{\Ov{\vt}} \Grad \mathfrak{T} \cdot \Grad \MTC      } \dt \br	
	&\leq -  \int_0^\tau \intO{  \Ov{\vr} \left(\frac{\partial s(\Ov{\vr}, \Ov{\vt})}{\partial \vr}(\mathfrak{R} - \mathcal{R}) + \frac{\partial s(\Ov{\vr}, \Ov{\vt})}{\partial \vt}(\mathfrak{T} - \MTC) \right) {\frac{\Ov{\vt} \alpha (\Ov{\vr}, \Ov{\vt} )}{ c_p (\Ov{\vr}, \Ov{\vt} )} } \Grad \Big( G + |\vc x_h|^2 \Big) \cdot \vU                       } \dt \br
	&\quad - { \int_0^\tau \intO{  \Ov{\vr} \left(\frac{\partial s(\Ov{\vr}, \Ov{\vt})}{\partial \vr}(\mathfrak{R} - \mathcal{R}) + \frac{\partial s(\Ov{\vr}, \Ov{\vt})}{\partial \vt}(\mathfrak{T} - \MTC) \right)  \frac{\kappa(\Ov{\vt})}{\Ov{\vr} c_p (\Ov{\vr}, \Ov{\vt} )} \Del \MTC                  } \dt} \br	
	&\quad { -  \int_0^\tau \intO{  \Ov{\vr} \left(\frac{\partial s(\Ov{\vr}, \Ov{\vt})}{\partial \vr}(\mathfrak{R} - \mathcal{R}) + \frac{\partial s(\Ov{\vr}, \Ov{\vt})}{\partial \vt}(\mathfrak{T} - \MTC) \right)  \frac{1}{\Ov{\vr} c_p (\Ov{\vr}, \Ov{\vt} )} \xi(t)                  } \dt} \br
	& \quad   +  \int_0^\tau \intO{ \lan \mathcal{R} - \mathfrak{R} \ran  \Grad \Big( G + \frac{1}{2}|\vc x_h|^2 \Big)  \cdot {\vU} } \dt  \br 
	&\quad + C \int_0^\tau \intO{E_\ep \left(\vre, \vte, \vue \Big| \Ov{\vr} + \ep \mathcal{R}, \Ov{\vt} + \ep \MTC, \vU \right) } \dt   + \omega_\ep(\tau).
	\label{A24}
\end{align}

\medskip

\noindent
$\bullet$ Now, as 
\[
\intO{ (\mathcal{R} - \mathfrak{R}) } = 0, 
\] 
we have 
\begin{align}
	-  \int_0^\tau \intO{  \Ov{\vr} \left(\frac{\partial s(\Ov{\vr}, \Ov{\vt})}{\partial \vr}(\mathfrak{R} - \mathcal{R}) + \frac{\partial s(\Ov{\vr}, \Ov{\vt})}{\partial \vt}(\mathfrak{T} - \MTC) \right)  \frac{1}{\Ov{\vr} c_p (\Ov{\vr}, \Ov{\vt} )} \xi(t)                  } \dt \nonumber \br 
	=- \int_0^\tau \intO{  \frac{\partial s(\Ov{\vr}, \Ov{\vt})}{\partial \vt} \left[
		\frac{\partial p(\Ov{\vr}, \Ov{\vt} )}{\partial \vr} \left( \frac{\partial p(\Ov{\vr}, \Ov{\vt} )}{\partial \vt} \right)^{-1}(\mathfrak{R} - \mathcal{R}) + (\mathfrak{T} - \MTC) \right]  \frac{1}{ c_p (\Ov{\vr}, \Ov{\vt} )} \xi(t)                  } \dt.
	\label{A25}
\end{align}
Note that both couples $\big(\mathcal{R} , \mathcal{T}\big)$ and $\big(\mathfrak{R}, \mathfrak{T}\big)$ satisfy the Boussinesq relation \eqref{TS3} and \eqref{A9}, respectively.
It follows in particular that the quantity 
\[
\left[
\frac{\partial p(\Ov{\vr}, \Ov{\vt} )}{\partial \vr} \left( \frac{\partial p(\Ov{\vr}, \Ov{\vt} )}{\partial \vt} \right)^{-1}(\mathfrak{R} - \mathcal{R}) + (\mathfrak{T} - \MTC) \right]
\]
is independent of $x$.

Similarly, we may rewrite 
\begin{align}
	-\int_0^\tau \intO{  \Ov{\vr} \left(\frac{\partial s(\Ov{\vr}, \Ov{\vt})}{\partial \vr}(\mathfrak{R} - \mathcal{R}) + \frac{\partial s(\Ov{\vr}, \Ov{\vt})}{\partial \vt}(\mathfrak{T} - \MTC) \right)  \frac{\kappa(\Ov{\vt})}{\Ov{\vr} c_p (\Ov{\vr}, \Ov{\vt} )} \Del \MTC                  } \dt \br 
	= -	\int_0^\tau \intO{  \frac{\partial s(\Ov{\vr}, \Ov{\vt})}{\partial \vr}\left[ (\mathfrak{R} - \mathcal{R}) +  \frac{\partial p(\Ov{\vr}, \Ov{\vt}) }{\partial \vt}
		\left(  \frac{\partial p(\Ov{\vr}, \Ov{\vt}) }{\partial \vr} \right)^{-1} (\mathfrak{T} - \MTC)  \right]  \frac{\kappa(\Ov{\vt})}{ c_p (\Ov{\vr}, \Ov{\vt} )} \Del \MTC                  } \dt \br 
	+ \int_0^\tau \intO{ \left(  \frac{\partial s(\Ov{\vr}, \Ov{\vt})}{\partial \vr}\frac{\partial p(\Ov{\vr}, \Ov{\vt}) }{\partial \vt}
		\left(  \frac{\partial p(\Ov{\vr}, \Ov{\vt}) }{\partial \vr} \right)^{-1} (\mathfrak{T} - \MTC) - \frac{\partial s(\Ov{\vr}, \Ov{\vt})}{\partial \vt} (\mathfrak{T} - \MTC)\right) \frac{\kappa(\Ov{\vt})} { c_p (\Ov{\vr}, \Ov{\vt} )} \Del \MTC  } \dt, 
	\label{A26}	
\end{align}
where, again, we have that
\[
\left[ (\mathfrak{R} - \mathcal{R}) +  \frac{\partial p(\Ov{\vr}, \Ov{\vt}) }{\partial \vt}
\left(  \frac{\partial p(\Ov{\vr}, \Ov{\vt}) }{\partial \vr} \right)^{-1} (\mathfrak{T} - \MTC)  \right]
\]
is independent of $x$.

Next, we may integrate equation \eqref{A23a} over the spatial domain $\Omega$, obtaining 
\begin{align} \left[ \left( \Ov{\vt} \alpha (\Ov{\vr} ,\Ov{\vt} ) 
	\frac{\partial p(\Ov{\vr}, \Ov{\vt})}{\partial \vt} \right)^{-1} -   \frac{1}{\Ov{\vr} c_p (\Ov{\vr}, \Ov{\vt}) } \right] |\Omega| \xi (t) = \int_{\partial \Omega} \frac{\kappa(\Ov{\vt})}{\Ov{\vr} c_p (\Ov{\vr}, \Ov{\vt} ) } \Grad \MTC \cdot \vc{n} \ \D \sigma_x.
	\label{A28}
\end{align}

Thus, substituting the integral in \eqref{A26} by \eqref{A28}, we can compute the sum of \eqref{A25} with the first integral in \eqref{A26}:
\begin{align}
	-  &\intO{  \frac{\partial s(\Ov{\vr}, \Ov{\vt})}{\partial \vt} \left[
		\frac{\partial p(\Ov{\vr}, \Ov{\vt} )}{\partial \vr} \left( \frac{\partial p(\Ov{\vr}, \Ov{\vt} )}{\partial \vt} \right)^{-1}(\mathfrak{R} - \mathcal{R}) + (\mathfrak{T} - \MTC) \right]  \frac{1}{ c_p (\Ov{\vr}, \Ov{\vt} )} \xi(t)                  } \br 
	-	 &\int_{\Omega}  \frac{\partial s(\Ov{\vr}, \Ov{\vt})}{\partial \vr}\left[ (\mathfrak{R} - \mathcal{R}) +  \frac{\partial p(\Ov{\vr}, \Ov{\vt}) }{\partial \vt}
	\left(  \frac{\partial p(\Ov{\vr}, \Ov{\vt}) }{\partial \vr} \right)^{-1} (\mathfrak{T} - \MTC)  \right] \times \br & \times \left[ \Ov{\vr} \left( \Ov{\vt} \alpha (\Ov{\vr} ,\Ov{\vt} ) 
	\frac{\partial p(\Ov{\vr}, \Ov{\vt})}{\partial \vt} \right)^{-1} -   \frac{1}{ c_p (\Ov{\vr}, \Ov{\vt}) } \right] \xi (t) \dx.	
	\label{S25}
\end{align}

Now, we use Gibbs' relation along with the specific formulae for $\alpha$ and $c_p$ to compute
\begin{align} 
	&- \frac{\partial s(\Ov{\vr}, \Ov{\vt})}{\partial \vt} 
	\frac{\partial p(\Ov{\vr}, \Ov{\vt} )}{\partial \vr} \left( \frac{\partial p(\Ov{\vr}, \Ov{\vt} )}{\partial \vt} \right)^{-1} - \frac{\partial s(\Ov{\vr}, \Ov{\vt} )}{\partial \vr} 
	\left( c_p (\Ov{\vr}, \Ov{\vt})\Ov{\vr} \left( \Ov{\vt} \alpha (\Ov{\vr} ,\Ov{\vt} ) 
	\frac{\partial p(\Ov{\vr}, \Ov{\vt})}{\partial \vt} \right)^{-1} - 1 \right) \br 
	= &- \frac{1}{\Ov{\vt}}\frac{\partial e(\Ov{\vr}, \Ov{\vt})}{\partial \vt} 
	\frac{\partial p(\Ov{\vr}, \Ov{\vt} )}{\partial \vr} \left( \frac{\partial p(\Ov{\vr}, \Ov{\vt} )}{\partial \vt} \right)^{-1} \br &- \frac{\partial s(\Ov{\vr}, \Ov{\vt} )}{\partial \vr}
	\left[ \Ov{\vr} \left( \frac{\partial e(\Ov{\vr}, \Ov{\vt}) }{\partial \vt} + \frac{1}{\Ov{\vr}} \Ov{\vt} \alpha (\Ov{\vr}, \Ov{\vt}) \frac{\partial p (\Ov{\vr}, \Ov{\vt} )}{\partial \vt} \right)\left( \Ov{\vt} \alpha (\Ov{\vr} ,\Ov{\vt} ) 
	\frac{\partial p(\Ov{\vr}, \Ov{\vt})}{\partial \vt} \right)^{-1}
	- 1 \right] \br 
	= &- \frac{1}{\Ov{\vt}}\frac{\partial e(\Ov{\vr}, \Ov{\vt})}{\partial \vt} 
	\frac{\partial p(\Ov{\vr}, \Ov{\vt} )}{\partial \vr} \left( \frac{\partial p(\Ov{\vr}, \Ov{\vt} )}{\partial \vt} \right)^{-1}  
	- \frac{\partial s(\Ov{\vr}, \Ov{\vt} )}{\partial \vr}
	\left[ \Ov{\vr} \frac{\partial e(\Ov{\vr}, \Ov{\vt}) }{\partial \vt}\left(
	\Ov{\vt} \alpha (\Ov{\vr} ,\Ov{\vt} ) 
	\frac{\partial p(\Ov{\vr}, \Ov{\vt})}{\partial \vt} \right)^{-1} 
	\right] \br 
	= &- \frac{1}{\Ov{\vt}}\frac{\partial e(\Ov{\vr}, \Ov{\vt})}{\partial \vt} 
	\frac{\partial p(\Ov{\vr}, \Ov{\vt} )}{\partial \vr} \left( \frac{\partial p(\Ov{\vr}, \Ov{\vt} )}{\partial \vt} \right)^{-1}  
	+ \frac{1}{\Ov{\vr}} \frac{\partial p(\Ov{\vr}, \Ov{\vt} )}{\partial \vt}
	\left[  \frac{\partial e(\Ov{\vr}, \Ov{\vt}) }{\partial \vt}\left(
	\Ov{\vt} \alpha (\Ov{\vr} ,\Ov{\vt} ) 
	\frac{\partial p(\Ov{\vr}, \Ov{\vt})}{\partial \vt} \right)^{-1} 
	\right] = 0.   
	\nonumber
\end{align}
We conclude that the coefficient multiplying the term $\mathcal{R} - \mathfrak{R}$ vanishes. In the same way, we may deduce that the coefficient 
multiplying $\mathcal{T} - \mathfrak{T}$ vanishes. 

Next, concerning the second term in \eqref{A26}, we observe that, using Gibbs' relation and our hypotheses on the constitutive relations, we can write
\[
\frac{\partial s (\Ov{\vr}, \Ov{\vt} ) }{\partial \vr} = - \frac{1}{\Ov{\vr}^2} \frac{\partial p(\Ov{\vr}, \Ov{\vt})} 
{\partial \vt};
\] 
whence 
\begin{align}  
	\left[ \frac{\partial s(\Ov{\vr}, \Ov{\vt})}{\partial \vr}\frac{\partial p(\Ov{\vr}, \Ov{\vt}) }{\partial \vt}
	\left(  \frac{\partial p(\Ov{\vr}, \Ov{\vt}) }{\partial \vr} \right)^{-1}  - \frac{\partial s(\Ov{\vr}, \Ov{\vt})}{\partial \vt} \right]  \frac{\kappa(\Ov{\vt})} { c_p (\Ov{\vr}, \Ov{\vt} )} \br =
	-\left[ \frac{1}{\Ov{\vr}^2} \left( \frac{\partial p(\Ov{\vr}, \Ov{\vt}) }{\partial \vt} \right)^2
	\left(  \frac{\partial p(\Ov{\vr}, \Ov{\vt}) }{\partial \vr} \right)^{-1}  + \frac{1}{\Ov{\vt}} \frac{\partial e(\Ov{\vr}, \Ov{\vt})}{\partial \vt} \right]  \frac{\kappa(\Ov{\vt})} { c_p (\Ov{\vr}, \Ov{\vt} )} = { -\frac{\kappa(\Ov{\vt})}{\Ov{\vt}}.}
	\nonumber
\end{align}

We conclude by collecting the integrals containing $\Grad \Big( G + \frac{1}{2}|\vc x_h|^2 \Big)$:
\begin{align}
	-  &\intO{  \Ov{\vr} \left(\frac{\partial s(\Ov{\vr}, \Ov{\vt})}{\partial \vr}(\mathfrak{R} - \mathcal{R}) + \frac{\partial s(\Ov{\vr}, \Ov{\vt})}{\partial \vt}(\mathfrak{T} - \MTC) \right) {\frac{\Ov{\vt} \alpha (\Ov{\vr}, \Ov{\vt} )}{ c_p (\Ov{\vr}, \Ov{\vt} )} } \Grad \Big(G + \frac{1}{2}|\vc x_h|^2 \Big) \cdot \vU                       }	\br 
	&+   \intO{ \lan \mathcal{R}  - \mathfrak{R} \ran \Grad \Big( G + \frac{1}{2}|\vc x_h|^2 \Big)  \cdot {\vU} } \br
	& =   -  \intO{  \Ov{\vr} \left(\frac{\partial s(\Ov{\vr}, \Ov{\vt})}{\partial \vr}\lan\mathfrak{R} - \mathcal{R}\ran + \frac{\partial s(\Ov{\vr}, \Ov{\vt})}{\partial \vt}\lan \mathfrak{T} - \MTC \ran \right) {\frac{\Ov{\vt} \alpha (\Ov{\vr}, \Ov{\vt} )}{ c_p (\Ov{\vr}, \Ov{\vt} )} } \Grad \Big(G + \frac{1}{2}|\vc x_h|^2 \Big) \cdot \vU                       }	\br
	&+   \intO{ \lan \mathcal{R}  - \mathfrak{R} \ran \Grad \Big( G + \frac{1}{2}|\vc x_h|^2 \Big)  \cdot {\vU} } \br 
	&=  \intO{  \left[ \Ov{\vr} \left(\frac{\partial s(\Ov{\vr}, \Ov{\vt})}{\partial \vr}\Grad  \lan\mathfrak{R} - \mathcal{R}\ran + \frac{\partial s(\Ov{\vr}, \Ov{\vt})}{\partial \vt}\Grad \lan\mathfrak{T} - \MTC\ran \right) {\frac{\Ov{\vt} \alpha (\Ov{\vr}, \Ov{\vt} )}{ c_p (\Ov{\vr}, \Ov{\vt} )} } \Big( G + \frac{1}{2}|\vc x_h|^2 \Big) \right] \cdot \vU                       } \br
	&-   \intO{ \Big( G + \frac{1}{2}|\vc x_h|^2 \Big) \Grad \lan \mathcal{R} - \mathfrak{R}\ran    \cdot {\vU} }.
	\nonumber
\end{align}
Furthermore, using Boussinesq relation,  we can compute
\begin{align}
	&\Ov{\vr} \left(\frac{\partial s(\Ov{\vr}, \Ov{\vt})}{\partial \vr}\Grad \lan \mathfrak{R} - \mathcal{R}\ran + \frac{\partial s(\Ov{\vr}, \Ov{\vt})}{\partial \vt}\Grad \lan \mathfrak{T} - \MTC \ran \right) {\frac{\Ov{\vt} \alpha (\Ov{\vr}, \Ov{\vt} )}{ c_p (\Ov{\vr}, \Ov{\vt} )} } \br 
	=&\Ov{\vr} \left(\frac{\partial s(\Ov{\vr}, \Ov{\vt})}{\partial \vr}\Grad \lan\mathfrak{R} - \mathcal{R}\ran - \frac{\partial s(\Ov{\vr}, \Ov{\vt})}{\partial \vt} { \frac{\partial p(\Ov{\vr}, \Ov{\vt}) }{\partial \vr} 
		\left( \frac{\partial p(\Ov{\vr}, \Ov{\vt}) }{\partial \vt} \right)^{-1} } \Grad \lan \mathfrak{R} - \mathcal{R} \ran \right) {\frac{\Ov{\vt} \alpha (\Ov{\vr}, \Ov{\vt} )}{ c_p (\Ov{\vr}, \Ov{\vt} )} }	\br 
	=& - \Ov{\vr} \left( \frac{1}{\Ov{\vr}^2} \frac{\partial p (\Ov{\vr}, \Ov{\vt}) }{\partial \vt} \Grad \lan\mathfrak{R} - \mathcal{R}\ran + \frac{1}{\Ov{\vt}} \frac{\partial e(\Ov{\vr}, \Ov{\vt})}{\partial \vt} { \frac{\partial p(\Ov{\vr}, \Ov{\vt}) }{\partial \vr} 
		\left( \frac{\partial p(\Ov{\vr}, \Ov{\vt}) }{\partial \vt} \right)^{-1} }  \Grad \lan \mathfrak{R} - \mathcal{R}\ran \right) {\frac{\Ov{\vt} \alpha (\Ov{\vr}, \Ov{\vt} )}{ c_p (\Ov{\vr}, \Ov{\vt} )} }\br 
	=& - \Grad \lan\mathfrak{R} - \mathcal{R}\ran.
	\nonumber
\end{align}

Summing up the previous relations and using the fact that $\mathcal{T}$ and $\mathfrak{T}$ share the same boundary values 
\[
\mathfrak{T}|_{\partial\Omega}=\MTC|_{\partial\Omega} = \vtB, 
\] 
we may rearrange terms in \eqref{A24} reaching the desired conclusion  
\begin{align}
	&\intO{ E_\ep \left(\vre, \vte, \vue \Big| \Ov{\vr} + \ep \mathcal{R}, \Ov{\vt} + \ep \MTC, \vU \right) (\tau, \cdot) }  \br 
	&+ \int_0^\tau \intO{ \Big( \mathbb{S} (\Ov{\vt}, \Ds \vue) - \mathbb{S} (\Ov{\vt}, \Ds \vU) \Big) : \Big( \Ds \vue - \Ds \vU  \Big) } \dt \br  &+ \int_0^\tau \intO{ \frac{\kappa (\Ov{\vt} ) }{\Ov{\vt}} \left|\Grad 
		\left(	\frac{\vte - \Ov{\vt}}{\ep} \right) - \Grad \MTC \right|^2 } \dt \br
	&\aleq \int_0^\tau \intO{E_\ep \left(\vre, \vte, \vue \Big|\Ov{\vr} + \ep \mathcal{R}, \Ov{\vt} + \ep \MTC, \vU \right) } \dt   + \omega_\ep(\tau).
	\label{A30}
\end{align}
Thus the standard application of the Gr\" onwall lemma yields the conclusion of Theorem \ref{TM1}.

\section{Concluding remarks} \label{s:comments}

\begin{enumerate}
	\item
The result holds also in the case when the underlying horizontal domain is not simply connected, in particular in the case of two concentric cylinders.
Indeed, this property was only used in \eqref{eq:potential} to deduce the existence of the potential $\Phi$ such that
\[
\Grad \Phi = - \vc e_3\times \vU. %(- u^2_h, u^1_h, 0). 
\]
However, as the field $(- u^2_h, u^1_h)$ vanishes on the boundary of the horizontal domain $B(r)$ it can be extended to be zero outside $B(r)$
and the corresponding potential $\Phi$ can be constructed.   

\item Imposing complete--slip boundary conditions at the horizontal boundaries of the domain $\Omega$ eliminates the effects of the
Ekman boundary layers, which do not appear in our context.
Replacing the complete--slip boundary conditions \eqref{B2} by Navier type boundary conditions 
\[
\vu \cdot \vc{n}|_{ B(r) \times \{ x_3 = 0, 1 \} } = 0,\ 
[\mathbb{S}(\vt, \Ds \vu) \cdot \vc{n} + \beta \vu ] \times \vc{n} |_{ B(r) \times \{ x_3 = 0, 1 \} } = 0
\]
would produce Ekman damping (associated with the so--called Ekman pumping phenomenon) in the limit momentum equation:
\[
\Ov{\vr} \Big[ \partial_t \vuh +  \Divh (\vuh \otimes \vuh) \Big] + \Gradh \Pi = \mu(\Ov{\vt}) \Delta_h \vuh - 2 \beta \vuh + \left< r \right> 
\Gradh\Big( G + \frac{1}{2}|\vc x_h|^2 \Big),
\] 
cf. Chemin et al \cite{CDGG}.
We leave to the interested reader to elaborate the necessary modifications in the proof.

\item For fast rotating fluids as the ones considered in the present work, the presence of vertical walls combined with the no-slip condition \eqref{B1}
typically entails the presence of boundary layers, also known as \emph{Munk boundary layers}. These are produced (see e.g. \cite[Chapter 11]{CDGG})
by a balance between the rotation and the pressure gradient at first order. Under the scaling considered in this paper, however, the singular perturbation
operator can be written (roughly, at least at first order) as
\[
\left\{ \begin{array}{l}
             \dfrac{1}{\sqrt{\ep}} \vc e_3\times \vu - \mu(\Ov{\vt}) \Del \vu + \dfrac{1}{\ep} \Grad r = 0 \\[2ex]
             \Div \vu = 0 .
        \end{array}
\right.
\]
Thus, the rotation is a lower order term and a balance between it and a suitable gradient term occurs only at higher order, which is however not captured by
the relative energy method. As a result, Munk boundary layers do not appear in our context.

\end{enumerate} 

\def\cprime{$'$} \def\ocirc#1{\ifmmode\setbox0=\hbox{$#1$}\dimen0=\ht0
	\advance\dimen0 by1pt\rlap{\hbox to\wd0{\hss\raise\dimen0
			\hbox{\hskip.2em$\scriptscriptstyle\circ$}\hss}}#1\else {\accent"17 #1}\fi}

%\bibliography{citace}
%\bibliographystyle{plain}

\end{document}